\documentclass[onefignum,onetabnum]{siamonline220329}

\usepackage{lipsum}
\usepackage{amsfonts}
    \usepackage{mathtools}
    \usepackage{mathrsfs}
    \usepackage{bbm}
\usepackage{graphicx}
\usepackage{caption}
\usepackage{subcaption}
\usepackage{tikz}

\usepackage{epstopdf}
\usepackage[]{algorithm}
\usepackage{algorithmic}
\ifpdf
  \DeclareGraphicsExtensions{.eps,.pdf,.png,.jpg}
\else
  \DeclareGraphicsExtensions{.eps}
\fi

\newsiamremark{remark}{Remark}
\newsiamremark{hypothesis}{Hypothesis}
\crefname{hypothesis}{Hypothesis}{Hypotheses}
\newsiamthm{claim}{Claim}

\headers{Long-Time Accuracy of the Ensemble Kalman Filter}{D. Sanz-Alonso and N. Waniorek}

\title{Long-Time Accuracy of Ensemble Kalman Filters \\ for Chaotic and Machine-Learned Dynamical Systems}

\author{Daniel Sanz-Alonso\thanks{Department of Statistics, University of Chicago, Chicago, IL 
  (\email{sanzalonso@uchicago.edu}).}
\and Nathan Waniorek\thanks{Committee on Computational and Applied Mathematics, University of Chicago, Chicago, IL 
  (\email{waniorek@uchicago.edu}).}
}

\usepackage{amsopn}

\providecommand{\mathbbm}{\mathbb} 

\renewcommand{\phi}{\varphi}
\newcommand{\eps}{\varepsilon}

\newcommand{\Tr}{{\rm{Tr}}}
\newcommand{\R}{\mathbbm{R}}
\newcommand{\E}{\mathbb{E}}

\newcommand{\N}{\mathbbm{N}}
\newcommand{\Z}{\mathbbm{Z}}

\renewcommand{\L}{\mathcal{L}}

\newcommand{\B}{\mathcal{B}}
\newcommand{\msK}{\mathscr{K}}

\definecolor{mygreen}{rgb}{0.13,0.55,0.13}

\newcommand{\Sam}{N}
\newcommand{\sam}{n}

\newcommand{\Nc}{\mathcal{N}}
\renewcommand{\H}{\mathcal{H}}

\usepackage[T1]{fontenc}
\usepackage[utf8]{inputenc}
\usepackage{graphicx}
\usepackage{placeins}
\usepackage{enumerate}
\usepackage{enumitem} 

\newtheorem{assumption}[theorem]{Assumption}
\newtheorem{example}[theorem]{Example}

\ifpdf
\hypersetup{
  pdftitle={Long-time accuracy of ensemble Kalman filters for chaotic dynamical systems and machine-learned surrogate models},
  pdfauthor={D. Sanz-Alonso and N. Waniorek}
}
\fi

\begin{document}

\maketitle

\begin{abstract}
Filtering is concerned with online estimation of the state of a dynamical system from partial and noisy observations. In applications where the state is high dimensional, ensemble Kalman filters are often the method of choice. This paper establishes long-time accuracy of ensemble Kalman filters. We introduce conditions on the dynamics and the observations under which the estimation error remains small in the long-time horizon. Our theory covers a wide class of partially-observed chaotic dynamical systems, which includes the Navier-Stokes equations and Lorenz models. In addition, we prove long-time accuracy of ensemble Kalman filters with surrogate dynamics, thus validating the use of machine-learned forecast models in ensemble data assimilation.
\end{abstract}



\section{Introduction}
The filtering problem involves the sequential estimation of the state of a dynamical system from partial and noisy observations. For high-dimensional problems, particularly those arising in the geophysical sciences \cite{kalnay2003atmospheric,oliver2008inverse,houtekamer1998data}, the ensemble Kalman filter \cite{evensen1994sequential} is often the method of choice. The algorithm propagates an ensemble of particles under the dynamics and nudges the ensemble to new observations using empirical covariances in a Kalman style update. 
In the linear-Gaussian setting, the ensemble Kalman filter converges to the Kalman filter as the ensemble size grows \cite{mandel2011convergence,kwiatkowski2015convergence}, thus fully characterizing the filtering distribution and yielding the optimal estimator of the state \cite{sanz2023inverse}. Furthermore, non-asymptotic theory shows that an ensemble size much smaller than the state dimension often suffices to approximate the Kalman filter \cite{al2024non}. For nonlinear dynamics, however, the ensemble Kalman filter does not reproduce the filtering distribution even in the large ensemble limit \cite{le2009large,ernst2015analysis,law2012evaluating}.  Despite this caveat, this paper rigorously justifies the use of ensemble Kalman filters in practical settings with nonlinear dynamics, partial observations, and an ensemble size smaller than the state dimension. We do so by providing conditions on the dynamics and the observations that guarantee the long-time accuracy of the ensemble Kalman filter as a state estimator. Our conditions cover a wide class of partially-observed chaotic dissipative dynamical systems, which includes the Navier-Stokes equations and Lorenz models. 

Motivated by recent advancements in machine learning emulation of geophysical systems, this paper additionally studies ensemble Kalman filters that use surrogate models for the dynamics. Examples of machine learning models for global weather forecasting include ClimaX \cite{nguyen2023climax}, FourCastNet \cite{pathak2022fourcastnet}, GraphCast \cite{lam2022graphcast}, Pangu \cite{bi2023accurate}, and FengWu \cite{chen2023fengwu}, among others. Due to the chaotic nature of the dynamical systems of interest, it is unreasonable to expect these surrogate models to provide accurate state predictions over long time scales, motivating work on learning surrogates that preserve the time-invariant statistical properties of the dynamics \cite{jiang2024training,park2024dynamical,tang2024learning}. However, in data assimilation, long-time filter accuracy may only need forecasts to be accurate over one assimilation cycle \cite{adrian2024data}. Consequently, cheap-to-evaluate machine learning surrogate models have the potential to enormously speed up data assimilation algorithms and allow for the generation of much larger ensembles than is possible with traditional numerical approximations of the true dynamics \cite{krasnopolsky2023using}. This paper establishes precise mathematical conditions for long-time accuracy of ensemble Kalman filters with surrogate models. 

\subsection{Related work}
\paragraph{Filter accuracy} 
A significant body of work on the accuracy of filtering algorithms builds on the theory of synchronization in dynamical systems \cite{pecora1990synchronization}. The seminal paper \cite{hayden2011discrete} establishes the convergence of discrete-time data assimilation methods for the Lorenz-63 and Navier-Stokes equations with partial but noiseless observations. These results were extended to more realistic observation models and noisy observations in \cite{azouani2014continuous,foias2016discrete,bessaih2015continuous}, and similar ideas are used to study cycled data assimilation schemes in \cite{moodey2013nonlinear}. Long-time accuracy for the 3DVar algorithm is shown in \cite{brett2013accuracy} for the Navier-Stokes equations, \cite{law2013analysis} for the Lorenz-63 system, and \cite{law2016filter} for the Lorenz-96 system. The accuracy of approximate Gaussian filters for the Navier-Stokes equations with model errors was studied numerically in \cite{branicki2018accuracy}, and conditions guaranteeing the long-time accuracy of 3DVar with surrogate models were given in \cite{adrian2024data}. A rigorous connection between the determining modes property of the Navier-Stokes equations and the synchronization of filters was made in \cite{carlson2024determining}, and the relationship between the nudging filter and synchronization filter for the Navier-Stokes equations was studied in \cite{carlson2024infinite}.

\paragraph{Ensemble Kalman filters} Much of the work on the theoretical foundations of ensemble Kalman filters has focused on the linear-Gaussian setting. In this setting, it is known that the ensemble Kalman filter converges to the Kalman filter in the large ensemble asymptotic \cite{mandel2011convergence,kwiatkowski2015convergence}. In the nonlinear setting, a similar line of work assumes the \textit{mean-field} perspective \cite{calvello2022ensemble,herty2019kinetic,ding2021ensemble,ding2021ensemblebis} and analyzes the large ensemble limits of the ensemble Kalman filter. Recent results using this approach have shown that the mean-field ensemble Kalman filter is near the filtering distribution in the near-Gaussian setting \cite{carrillo2024mean} and in the near-linear setting \cite{calvello2024accuracy}. These results only hold over a finite-time horizon and do not guarantee accurate state estimation. A complementary line of work provides a non-asymptotic analysis of ensemble Kalman updates, explaining how the ensemble Kalman filter can well approximate the mean-field updates with an ensemble size much smaller than the dimension of the state space \cite{al2024non,sanz2024analysis,al2024ensemble}. This perspective helps explain why localization techniques are so successful within ensemble Kalman filters and extends notions of localization and sparse covariance estimation to the function space setting \cite{al2023covariance,al2024optimal}. However, none of these non-asymptotic analyses hold over an infinite-time horizon, and they do not consider state estimation accuracy. A number of previous works have considered the state estimation accuracy of the ensemble Kalman filter, largely using the theory of synchronization in dynamical systems \cite{pecora1990synchronization}. The work \cite{kelly2014well} proves stability and accuracy of the ensemble Kalman filter with covariance inflation, but only in the setting with complete observations of the system. They provide numerical evidence of filter accuracy with partial observations, but do not extend the theorems to the partially-observed setting. Recently, \cite{takeda2024uniform} proves accuracy of the ensemble Kalman filter with multiplicative covariance inflation, but again requires full observations and the results are restricted to finite-dimensional dynamical systems. Notable works in the partially-observed setting include \cite{tong2015nonlinear,gonzalez2013ensemble,biswas2024unified}. In \cite{tong2015nonlinear}, the authors prove long-time stability of the ensemble Kalman filter with an adaptive variance inflation scheme, but do not prove filter accuracy. For hyperbolic dynamical systems, \cite{gonzalez2013ensemble} proves accuracy of the ensemble Kalman filter with partial observations using a shadowing argument. The work most similar to ours is \cite{biswas2024unified}, which proves accuracy of a continuous-time variant of the ensemble Kalman filter for the Navier-Stokes equations in the partially-observed setting. Their algorithm includes variance inflation and a somewhat atypical localization procedure where the unobserved components of the forecast covariance are projected to be zero. In contrast, our results hold in the discrete-time setting and do not require this localization procedure.

\paragraph{Machine learning surrogate models in data assimilation} Machine learning methods have made far-reaching impacts on data assimilation; see 
\cite{bocquet2023surrogate,cheng2023machine,bach2024inverse} and \cite[Chapter 10]{chen2023stochastic} for recent reviews. The use of machine learning surrogates has the potential to rapidly accelerate data assimilation algorithms, due to both the orders of magnitude speedup in simulating the dynamics and the possibility of employing automatic differentiation for the computation of adjoint models \cite{hatfield2021building,maulik2022efficient,xiao2023fengwu}, although these adjoint models may display unphysical behavior \cite{tian2024exploringusemachinelearning}. Within ensemble data assimilation schemes, surrogate dynamics models make it possible to generate much larger ensembles than traditional expensive numerical solvers  \cite{chattopadhyay2022towards,gottwald2021supervised} and can in turn lead to improved estimates of the background covariance \cite{chattopadhyay2023deep}. In the context of global weather prediction, case studies in \cite{adrian2024data} and \cite{kotsuki2024integrating} provide compelling numerical evidence that entirely data-driven surrogates can be used successfully within 3DVar and the ensemble Kalman filter, respectively. For the 3DVar algorithm, \cite{adrian2024data} provides conditions on the surrogate model that ensure accurate state estimation over long time horizons. An alternative approach is to use machine learning models to learn and correct the model error incurred by a traditional physics-based model \cite{chen2021bamcafe,farchi2021using}. Other approaches learn the state trajectory and model/model error together \cite{bocquet2020bayesian,brajard2020combining,chen2022autodifferentiable,gottwald2021combining,hamilton2016ensemble,nguyen2019like}. Another line of work uses machine learning models not as surrogates for the dynamics, but to learn other components of the data assimilation process such as the forecast error covariance, the analysis update, or the entire data assimilation procedure \cite{arcucci2021deep,mou2023combining,penny2022integrating,tsuyuki2022nonlinear}.

\subsection{Main contributions and outine}
\begin{itemize}
    \item Section \ref{sec:mainresults} formalizes the problem setting and states our main results. Our first main result, Theorem \ref{th:main1}, proves the long-time accuracy of the square-root ensemble Kalman filter with appropriate variance inflation given that the observations are sufficiently informative. The second main result, Theorem \ref{th:main2}, proves the long-time accuracy of the square-root ensemble Kalman filter with a surrogate model, provided the model is sufficiently accurate in the unobserved part of the system.
    \item Section \ref{sec:proofmain1} contains the proof of Theorem \ref{th:main1}. We first establish accuracy for a mean-field version of the ensemble Kalman filter in Theorem \ref{theorem:general mean field filter accuracy}. Then, we show that the ensemble means from the square-root ensemble Kalman filter remain close to those output by the idealized mean-field version in Theorem \ref{theorem:general finite ensemble filter accuracy}. 
    \item Section \ref{sec:proofmain2} contains the proof of Theorem \ref{th:main2}. As in Section \ref{sec:proofmain1}, we will analyze an idealized mean-field algorithm as an intermediate step in the analysis.
    \item Section \ref{sec:numerics} contains numerical experiments illustrating our theoretical results for the Lorenz-96 system. In particular, we show that long-time filter accuracy is achievable with machine-learned surrogate models that only provide accurate short-term forecasts.
    \item Section \ref{sec:conclusions} concludes with a summary of our work and presents several open directions for future research.
\end{itemize}

\section{Main results}\label{sec:mainresults}

\subsection{Long-time accuracy of ensemble Kalman filters}
In this subsection, we show long-time accuracy for a standard implementation of the ensemble Kalman filter. Our setting covers a wide range of partially-observed dissipative chaotic dynamical systems often used to benchmark filtering algorithms.  
\subsubsection{Set-up}\label{ssec:setting}
We study the filtering problem of sequentially estimating a time-evolving hidden signal as new observations become available.  
Suppose that the signal process $\{u_j\}_{j=1}^\infty \subset \H$ lies on a separable Hilbert space $(\mathcal{H},\langle\cdot,\cdot\rangle, \|\cdot \|)$  and its evolution is governed by a discrete-time dynamical system given by $u_0 \in \mathcal{H}$ and
\begin{align}\label{eq:dynamics}
\begin{split}
    u_{j} &=\Psi(u_{j-1}), \qquad j = 1, 2, \ldots
\end{split}    
\end{align}
The map $\Psi: \H \to \H$ is assumed to be known, but the initial condition $u_0 \in \H$ ---and hence the signal process--- are assumed to be unknown. 

To estimate the signal, we are given observations $\{y_j\}_{j=1}^\infty \subset \R^k$ that comprise finite-dimensional noisy measurements of the state 
\begin{align*}
    y_{j}=Hu_{j}+\eps \eta_j, \qquad j = 1, 2, \ldots
\end{align*}
Here, $H:\H \to \R^k$ is assumed to satisfy $H H^* = I_k$ and we define $P:=H^* H:\H\to \H$, which is an orthogonal projection onto a $k$-dimensional subspace of $\H$. We are interested in the case where $k \ll \text{dim}(\H) \in \N \cup \{\infty\}$
and assume that the observation errors $\{\eta_j\}_{j=1}^\infty$ are independent, with $\mathbb{E}[\eta_j]=0$ and $\text{Cov}[\eta_j]=R\in \R^{k\times k}.$ The parameter $\eps>0$ controls the size of the observation noise, which is often small in applications. Filtering concerns sequentially estimating the state $u_j$ given all the observations $\{y_i\}_{i=1}^j$ available at time $j.$

The ability of filtering algorithms to accurately estimate the signal over a long time window necessarily hinges on control-theoretic conditions on the signal and observation processes. We will adopt the following assumption, also used in \cite{sanz2015long}, on the dynamics and the observations. As discussed in Subsection \ref{ssec:examples} below, this assumption is satisfied by important test filtering problems, including partially-observed Lorenz models and Navier-Stokes equations. For a constant $\beta\geq 0,$ we define $V:\H\to [0,\infty)$ by
        \begin{align*}
            V(u):=\Bigl(\|u\|^2+\beta \|Pu\|^2\Bigr)^{1/2}.
        \end{align*}
The norm $V$ will act as a Lyapunov function in the theory.        
\begin{assumption}\label{assumption:ball and squeezing}
The following properties hold: \\
        {\sc (Absorbing ball property)} There exist  a constant $r>0$ such that the ball
        \begin{align}
            \B:=\left\{u\in \H : \|u\|\leq r\right\}
        \end{align}
         is absorbing and forward invariant for the dynamical system \eqref{eq:dynamics}.  \\
                 {\sc (Local Lipschitz continuity)}
        There exists $L>0$ such that, for any $u,v\in \B,$
        \begin{align}
            V^2\bigl(\Psi(u)-\Psi(v) \bigr)\leq LV^2(u-v).
        \end{align}
         { \sc (Squeezing property) }
        There exists $\alpha\in (0,1)$ such that, for any $u,v\in \B,$
        \begin{align}\label{eq:squeezingproperty}
            V^2\Bigl( (I-P) \bigl(\Psi(u)-\Psi(v)\bigr) \Bigr)\leq \alpha V^2\left(u-v \right).
        \end{align}
\end{assumption}

The absorbing ball property and local Lipschitz continuity are naturally satisfied by dissipative dynamical systems where there is an energy-loss mechanism. In such systems, the bounded invariant sets of the dynamics contain all the information about asymptotic behavior \cite[Section 1.8]{stuart1998dynamical}; therefore, these sets play a key role in the study of long-time filter accuracy. The squeezing property can be interpreted as a detectability condition that ensures that the dynamics contract in the unobserved part of the system \cite{sanz2015long}.

\subsubsection{Square-root ensemble Kalman filter}
Ensemble Kalman filters sequentially estimate the state using an ensemble of interacting particles. At time $j = 0,$ we generate an \emph{initial ensemble} $\mathsf{U}_{0} = \{ u_{0}^{(n)} \}_{n=1}^N$ of $N$ particles; for concreteness, we assume that $u_0^{(n)} \stackrel{\text{i.i.d.}}\sim \Nc(m_0,C_0), $ $1 \le n \le N,$ for given mean $m_0 \in \H$ and covariance operator  $C_0.$  Then, at each time $j \ge 1,$ we estimate the hidden state $u_j$ by the empirical mean $\widehat{m}_j$ of an ensemble  $\mathsf{U}_{j} = \{ u_{j}^{(n)} \}_{n=1}^N$. Particles are updated from time $j-1$ to time $j$ in two steps, referred to as \emph{prediction} and \emph{analysis}, which utilize respectively the dynamical system \eqref{eq:dynamics} and the observation $y_j$ acquired at time $j.$ Thus, each update can be schematically represented by:

\begin{equation}
    \mathsf{U}_{j-1} = \bigl\{ u_{j-1}^{(n)} \bigr\}_{n=1}^N \xrightarrow[\text{Dynamical System}]{\text{Prediction}} 
    \mathsf{V}_{j} = \bigl\{ v_{j}^{(n)} \bigr\}_{n=1}^N \xrightarrow[\text{Observation}]{\text{Analysis}} 
    \mathsf{U}_{j} = \bigl\{ u_{j}^{(n)} \bigr\}_{n=1}^N\, .
\end{equation}

In the prediction step, we use the dynamics governing the evolution of the state and set
\begin{equation}\label{eq:predictionstep}
    v_{j}^{(\sam)}= \Psi \bigl(P_{\B_V}u_{j-1}^{(\sam)}\bigr)+ \xi_{j}^{(n)}, \qquad \xi_{j}^{(n)} \stackrel{\text{i.i.d.}}{\sim}  \Nc(0,Q), \qquad 1 \le n \le N.
\end{equation}
Here, we follow \cite{sanz2015long} and project each particle $u_{j-1}^{(n)}$ into the absorbing ball $\mathcal{B}$ in Assumption \ref{assumption:ball and squeezing}; specifically, $P_{\mathcal{B}_V}$ denotes the orthogonal projection onto $\mathcal{B}$ in the $V$-norm. 
To promote spread of the prediction ensemble and avoid downward bias of estimates of the covariance of the state \cite{furrer2007estimation}, we introduce additional randomness $\xi_{j}^{(n)}$ in the propagation of each particle $u_{j-1}^{(n)}$. This is a standard procedure in data assimilation known as \emph{covariance inflation}; see for instance \cite[Chapter 10]{evensen2022data}, \cite[Section 4]{katzfuss2016understanding}, and the discussion in \cite{bach2024inverse}. 
We take the covariance inflation operators to be of the form $Q=aP$ for some constant $a>0.$ Since $P$ is an orthogonal projection onto a $k$-dimensional subspace, $Q$ is trace-class and positive semi-definite, and hence a valid covariance operator in $\H$. It will be relevant to note that for $\xi\sim \Nc(0,Q)$, it holds that  $\text{Cov}[H\xi]=HQH^*=aHPH^*=aI_{k\times k}.$

In the analysis step, we update the prediction ensemble $\mathsf{V}_{j} = \{ v_{j}^{(n)} \}_{n=1}^N$ into an analysis ensemble $\mathsf{U}_{j} = \{ u_{j}^{(n)}\}_{n=1}^N$ by imposing that the empirical mean and covariances are consistent with Kalman-type formulae for linear-Gaussian inverse problems. To that end, we define the \emph{Kalman gain} associated to covariance operator $\Sigma$ by
\begin{equation}\label{eq:Kalmangain}
    \msK(\Sigma):= \Sigma H^*(H\Sigma H^* + \eps R)^{-1}.
\end{equation}
Then, denoting 
by $(\widehat{\mu}_j, \widehat{\Sigma}_j)$ the empirical mean and covariance of $\mathsf{V}_{j}$  and
by $(\widehat{m}_j, \widehat{C}_j)$ the empirical mean and covariance of $\mathsf{U}_{j},$ we transform the particles so that
\begin{align}\label{eq:analysisstep}
\begin{split}
    \widehat{m}_{j}&=\widehat{\mu}_{j}+
        \msK(\widehat{\Sigma}_{j})\left(y_{j+1} -H\widehat{\mu}_{j}\right),\\
        \widehat{C}_{j} &
        = \bigl(I - \msK(\widehat{\Sigma}_{j}) H \bigr) \widehat{\Sigma}_{j}.  
        \end{split}
\end{align}
As described in \cite{tippett2003ensemble} and reviewed in \cite{al2024non}, many deterministic maps can be used to transform the prediction ensemble $\mathsf{V}_{j}$ into an analysis ensemble $\mathsf{U}_j$ that satisfies the moment constraints in \eqref{eq:analysisstep}, leading to different square-root ensemble Kalman filters. Notable examples include the Ensemble Transform Kalman Filter \cite{bishop2001adaptive} and the Ensemble Adjustment Kalman Filter \cite{anderson2001ensemble}. Our theory is agnostic to the choice of transformation, and holds for the general formulation of square-root ensemble Kalman filters summarized in Algorithm \ref{algorithm:finite ensemble}.

\begin{algorithm}[H]
\caption{Square-Root Ensemble Kalman Filter \label{algorithm:finite ensemble}}
\begin{algorithmic}[1]
\STATE {\bf Initialization}: Draw
\vspace{-0.1cm}
$$ u_0^{(n)} \stackrel{\text{i.i.d.}}{\sim} \Nc(m_0,C_0), \qquad  1 \le n \le N.$$
\vspace{-0.5cm}
\STATE For $j = 1, 2, \ldots$ do the following:
\STATE {\bf \index{prediction}Prediction}: Compute
\vspace{-0.3cm}
\begin{align*}
\begin{split}
v_{j}^{(\sam)}&= \Psi \bigl(P_{\B_V}u_{j-1}^{(\sam)} \bigr)+ \xi_{j}^{(n)}, \qquad \xi_{j}^{(n)} \stackrel{\text{i.i.d.}}{\sim}  \Nc(0,Q), \qquad 1 \le n \le N,  \\
\widehat{\mu}_{j} &= \frac{1}{\Sam}\sum^{\Sam}_{\sam=1} v_{j}^{(\sam)},  \qquad 
\widehat{\Sigma}_{j} = \frac{1}{\Sam}\sum^{\Sam}_{\sam=1}\bigl(v^{(\sam)}_{j}-\widehat{\mu}_{j}\bigr)\otimes \bigl(v^{(\sam)}_{j}-\widehat{\mu}_{j}\bigr).
\end{split}
\end{align*}
\vspace{-0.4cm}
\STATE{{{\bf \index{analysis}Analysis}}}: Transform $\{v_{j}^{(n)}\}_{n=1}^N\mapsto \{u_{j}^{(n)}\}_{n=1}^N$ such that
\vspace{-0.3cm}
	\begin{align*}
	\begin{split}
        \widehat{m}_{j}&:= \frac{1}{N}\sum^{N}_{n=1} u_{j}^{(n)}=\widehat{\mu}_{j}+
        \msK(\widehat{\Sigma}_{j})\left(y_{j+1} -H\widehat{\mu}_{j}\right),\\
        \widehat{C}_{j} &:= \frac{1}{N-1}\sum_{n=1}^N \bigl(u^{(n)}_{j}-\widehat{m}_{j}\bigr)\otimes \bigl(u^{(n)}_{j}-\widehat{m}_{j}\bigr) 
        = \bigl(I - \msK(\widehat{\Sigma}_{j}) H \bigr) \widehat{\Sigma}_{j}.
\end{split}
\end{align*}
\vspace{-0.4cm}
\STATE {\bf Output}: Means and covariances  $\{\widehat{m}_j, \widehat{C}_j\}_{j=1}^\infty.$
\end{algorithmic}
\end{algorithm}

We are ready to state our first main result, which establishes long-time accuracy of the ensemble Kalman filter. We show that both observed and unobserved parts of the state can be estimated at the level $\eps$ of the observation noise. The proof will be deferred to Section \ref{sec:proofmain1}.
\begin{theorem}\label{th:main1}
    Suppose that Assumption \ref{assumption:ball and squeezing} holds and that $N \ge 6k.$ Then, if the inflation parameter $a>0$ is sufficiently large, there is a constant $\mathsf{C}$ independent of $\eps$ such that
    \begin{align*}
        \lim_{j\to \infty}\sup \E \|\widehat{m}_j-u_j\|\leq \mathsf{C}\eps.
    \end{align*}
\end{theorem}

\begin{remark}
   Our theory requires the ensemble size $N$ to scale with the dimension $k$ of the observations, but \emph{not} with the state dimension. The assumption that $N \ge 6k$ is a requirement of the results in \cite{mourtada2022exact} that we use to lower bound the smallest eigenvalue of $H\widehat{\Sigma}_jH^*$.
   Instead of performing variance inflation by randomly perturbing each particle in \eqref{eq:predictionstep}, one could consider deterministic variance inflation by noiselessly propagating each particle through the dynamics, i.e. setting $v_{j}^{(\sam)}= \Psi \bigl(P_{\B_V}u_{j-1}^{(\sam)} \bigr) ,$ and then deterministically and additively inflating the empirical covariance. 
  Then, our proof techniques suggest that $N\geq k$ suffices. 
\end{remark}

\subsubsection{Application to partially-observed dissipative chaotic dynamical systems}\label{ssec:examples}
Here we verify Assumption \ref{assumption:ball and squeezing} for partially-observed signals arising from dissipative systems with quadratic, energy-conserving nonlinearity of the form
\begin{equation}\label{eq:dissipativesystem}
    \frac{du}{dt} + \mathsf{A} u + \mathsf{B}(u,u) = \mathsf{F}, \quad u(0) = u_0 \in \mathcal{H}.
\end{equation}
We denote by $\Delta t >0$ the time between observations and let $\Psi$ denote the $\Delta t$-flow of \eqref{eq:dissipativesystem}; that is, $\Psi(u_0)$ is the value at time $\Delta t$ of the solution to \eqref{eq:dissipativesystem} with initial condition $u_0 \in \mathcal{H}.$ With this set-up, the state $u_j$ in \eqref{eq:dynamics} represents the solution to \eqref{eq:dissipativesystem} at time $t_j:= j \Delta t$ with initial condition $u_0 \in \mathcal{H}.$ We next verify Assumption \ref{assumption:ball and squeezing} for partially-observed Lorenz-63, Lorenz-96, and Navier-Stokes equations, all of which can be written in the form \eqref{eq:dissipativesystem}. These examples are not new \cite{hayden2011discrete,robinson2001infinite,law2016filter,sanz2015long}, but we review them here to demonstrate the wide applicability of our novel long-time accuracy theory for ensemble Kalman filters.

\begin{example}\label{ex:Lorenz63}
    The Lorenz-63 model\footnote{In these examples, subscripts denote coordinates and coefficient indices rather than time.} with standard parameter values $(10,8/3,28)$ can be written in the form \eqref{eq:dissipativesystem} with $\H := \R^3,$ 
\begin{align}
    \mathsf{A} := \begin{bmatrix}
10 & -10 & 0\\
10 & 1 & 0 \\
0 & 0 & 8/3 
\end{bmatrix}, \qquad 
\mathsf{B}(u, \tilde{u}) :=
\begin{bmatrix}
 0\\
(u_1 \tilde{u}_3 + u_3 \tilde{u}_1 )/2 \\
-(u_1 \tilde{u}_2 + u_2 \tilde{u}_1 )/2
\end{bmatrix},
\qquad 
\mathsf{F} =
\begin{bmatrix}
 0\\
0 \\
-304/3
\end{bmatrix}.
\end{align}
The Lorenz-63 model satisfies the absorbing ball property, see e.g. \cite[Example 2.6]{law2015data}. Moreover, the squeezing property is satisfied observing only the first component of the state: for sufficiently small $\Delta t,$ the $\Delta t$-flow  $\Psi: \R^3 \to \R^3$ of the Lorenz equation and the observation map $H = \begin{bmatrix}
    1 & 0 & 0 
\end{bmatrix}:  \R^3 \to \R$ satisfy \eqref{eq:squeezingproperty} with $\beta=1$. We refer to \cite{sanz2015long} for a proof. 
\end{example}

\begin{example}\label{ex:Lorenz96}
     The Lorenz-96 model can be written in the form \eqref{eq:dissipativesystem} with $\H:=\R^d$ where we assume $d\in 3\N$,
     \begin{align}
         \mathsf{A} := I_{d\times d}, \qquad \mathsf{B}(u,\tilde{u}) := -\frac{1}{2}\begin{bmatrix}
             \tilde{u}_2u_d+u_2\tilde{u}_d-\tilde{u}_du_{d-1}-u_d\tilde{u}_{d-1}\\
             \vdots \\
             \tilde{u}_{i-1}u_{i+1}+u_{i-1}\tilde{u}_{i+1}-\tilde{u}_{i-2}u_{i-1}-u_{i-2}\tilde{u}_{i-1}\\
             \vdots \\
             \tilde{u}_{d-1}u_1+u_{d-1}\tilde{u}_1-\tilde{u}_{d-2}u_{d-1}-u_{d-2}\tilde{u}_{d-1}
         \end{bmatrix}, \qquad \mathsf{F} =
\begin{bmatrix}
 8\\
\vdots \\
8
\end{bmatrix}.
     \end{align}
     The Lorenz-96 model satisfies the absorbing ball property, see e.g. \cite[Example 2.7]{law2015data}. Moreover, it satisfies the squeezing property observing only two out of every three components of the state. Let  $H\in \R^{\frac{2}{3}d\times d}$ denote the identity matrix with every third row removed. Then,  for sufficiently small $\Delta t$, the $\Delta t$-flow $\Psi:\R^d\to \R^d$ of the Lorenz-96 model and the observation map $H$ satisfy \eqref{eq:squeezingproperty} with $\beta=1$. We refer to \cite{sanz2015long} for a proof.
\end{example}

\begin{example}\label{ex:NSE}
       We consider the incompressible Navier-Stokes equations on the 2-dimensional torus, $\mathbb{T}^2=[0,1)\times [0,1)$. Let $\mathcal{U}_0$ be the space of zero-mean, divergence-free, vector-valued polynomials $u$ from $\mathbb{T}^2$ to $\R^2$, and let $\mathcal{U}$ be the closure of $\mathcal{U}_0$ with respect to the $L^2$ norm. We define $P_{\mathcal{U}}: \bigl(L^2(\mathbb{T}^2)\bigr)^2 \to \mathcal{U}$ to be the Leray-Helmholtz orthogonal projector \cite{temam1995navier}. Functions in $\mathcal{U}$ can be written as a Fourier series
       \begin{align}
           u=\sum_{\ell\in \mathcal{L}}u_\ell e^{i\ell x}, \qquad \mathcal{L}=\Bigl\{2\pi(n_1,n_2):n_i\in \Z, (n_1,n_2)\neq (0,0)\Bigr\}.
       \end{align}
       We define $\H$ to be the closure of $\mathcal{U}_0$ with respect to norm defined by 
       \begin{align*}
           \|u\|^2:=\sum_{\ell\in \mathcal{L}} |\ell|^2|u_\ell|^2,
       \end{align*}
       which is equivalent to the Sobolev $H^1$ norm.
       The Navier-Stokes equations can be written in the form \eqref{eq:dissipativesystem} with
       \begin{align}
           \mathsf{A}:=-\nu P_{\mathcal{U}} \Delta, \qquad \mathsf{B}(u,\tilde{u})=P_{\mathcal{U}}(u \cdot \nabla \tilde{u}), \qquad \mathsf{F}=f\in \mathcal{U}.
       \end{align}
       Given an initial condition in $\H$, the Navier-Stokes equations have a unique strong solution for any time $\Delta t>0$ that depends continuously on the initial data in the $\H$ norm \cite[Theorem 9.5]{robinson2001infinite}. One can show that there exists $\theta=\theta(\nu)>0$ such that for every $u\in \H$, $\|\mathsf{A}u\|^2\geq \theta \|u\|^2$ \cite[Chapter 12]{robinson2001infinite}. Then, the absorbing ball property holds with 
       \begin{align*}
           \B=\biggl\{u\in \H : \|u\| \leq r:=\sqrt{2}\frac{\|f\|}{\theta}\biggr\}.
       \end{align*}
       We define $\mathcal{\L}_\lambda = \{  \ell \in \mathcal{L}: |\ell|^2 \le \lambda \}$ and set $k = \rm{card}(\mathcal{\L}_\lambda).$ Consider the observation operator 
            $H_{\lambda}u := (|\ell| u_\ell)_{\ell \in \mathcal{L}_\lambda}$
       with corresponding operator $P_{\lambda}:\H\to \H$ given by
       \begin{align*}
           P_{\lambda}u=\sum_{\ell \in \mathcal{L}_\lambda} u_\ell e^{i\ell x},
       \end{align*}
       which is an orthogonal projection onto a $k$-dimensional subspace of $\mathcal{H}.$
       There exist constants $\alpha=\alpha(r)\in (0,1)$ and $\lambda_*=\lambda_*(r)>0$ such that for all $\lambda>\lambda_*$, there exists $\Delta t_*=\Delta t_*(r,\lambda)$ such that, for all $u,v\in \B$ and all $\Delta t\leq \Delta t_*,$
       \begin{align*}
           \bigl\|(I-P) \Psi(u)-\Psi(v)\bigr\|\leq \alpha \bigl\|u-v\bigr\|.
       \end{align*}
       Hence, the squeezing property is satisfied with $\beta=0$ by the observation map $H=H_{\lambda}$ for any $\lambda\geq \lambda_*.$ We refer to \cite{sanz2015long} for a proof of the result as stated here, and to \cite[Chapter 14]{robinson2001infinite} for additional background on the squeezing property for the Navier-Stokes equations. 
       
\end{example}

\subsection{Long-time accuracy with surrogate models}
Next, we analyze ensemble Kalman filters that utilize a surrogate model for the dynamics.
\subsubsection{Set-up}
 We assume that the signal and observation processes are defined as in Subsection \ref{ssec:setting}, with dynamics map $\Psi$ and observation matrix $H$ satisfying Assumption \ref{assumption:ball and squeezing}.  We will study ensemble Kalman filters that utilize a surrogate dynamics model $\Psi^s$ satisfying:
\begin{assumption}\label{assumption:surrogate model assumption}
The following properties hold: \\
    {\sc (Bounded approximation error on $\mathcal{B}$)} There exists $\kappa>0$ such that, for any $u\in \mathcal{B},$
        \begin{align*}
        \|\Psi^s(u)-\Psi(u)\|\leq \kappa.
    \end{align*}
    {\sc (Local Lipschitz continuity)} There exists $L_s>0$ such that, for any $u,v\in \B,$
    \begin{align*}
        V^2 \bigl(\Psi^s(u)-\Psi^s(v)\bigr)\leq L_s  V^2  \bigl(u-v\bigr).
    \end{align*}
    {\sc(Accuracy in unobserved part)} There exists small $\delta>0$ such that
        \begin{align*}
        \bigl\|(I-P)\bigl(\Psi^s(u)-\Psi(u) \bigr) \bigr\|\leq \delta.
    \end{align*}
\end{assumption}
    
Assumption \ref{assumption:surrogate model assumption} controls the error of the surrogate model in the unobserved part of the system. Notice that the surrogate  model is only assumed to be accurate over a single discrete time-step, without requiring stability over a long time window. Consequently, our theory will validate the use of ensemble data assimilation with machine-learned surrogate models that provide accurate \emph{short-term} forecasts in the unobserved part of the system. 

\subsubsection{Square-root ensemble Kalman filter with a surrogate model}
We will study the following algorithm, obtained by replacing in Algorithm \ref{algorithm:finite ensemble} the true dynamics map $\Psi$ with the surrogate model $\Psi^s$:
\begin{algorithm}[H]
\caption{Square-Root Ensemble Kalman Filter with a Surrogate Model\label{algorithm:finite ensemble surrogate}}
\begin{algorithmic}[1]
\STATE {\bf Initialization}: Draw
\vspace{-0.1cm}
$$ u_0^{(n),s} \stackrel{\text{i.i.d.}}{\sim} \Nc(m_0,C_0), \qquad  1 \le n \le N.$$
\vspace{-0.5cm}
\STATE For $j = 1, 2, \ldots $ do the following:
\STATE {\bf \index{prediction}Prediction}: Compute
\vspace{-0.3cm}
\begin{align*}
\begin{split}
v_{j}^{(\sam),s}&= \Psi^s \bigl(P_{\B_V}u_{j-1}^{(\sam),s} \bigr)+ \xi_{j}^{(n)},  \qquad \xi_{j}^{(n)} \stackrel{\text{i.i.d.}}{\sim}  \Nc(0,Q), \qquad 1 \le n \le N, \\
\widehat{\mu}_{j}^s &= \frac{1}{\Sam}\sum^{\Sam}_{\sam=1} v_{j}^{(\sam),s},  \qquad 
\widehat{\Sigma}_{j}^s = \frac{1}{\Sam}\sum^{\Sam}_{\sam=1}\bigl(v^{(\sam),s}_{j}-\widehat{\mu}_{j}^s\bigr)\otimes \bigl(v^{(\sam),s}_{j}-\widehat{\mu}_{j}^s\bigr).
\end{split}
\end{align*}
\vspace{-0.35cm}
\STATE{{{\bf \index{analysis}Analysis}}}: Transform $\{v_{j}^{(n),s}\}_{n=1}^N\mapsto \{u_{j}^{(n),s}\}_{n=1}^N$ such that
\vspace{-0.3cm}
	\begin{align*}
	\begin{split}
        \widehat{m}^s_{j}&:= \frac{1}{N}\sum^{N}_{n=1} u_{j}^{(n),s} =\widehat{\mu}_{j}^s+
        \msK(\widehat{\Sigma}_{j}^s) 
         \left(y_{j} -H\widehat{\mu}^s_{j}\right),\\
        \widehat{C}^s_{j} &:= \frac{1}{N-1}\sum_{n=1}^N \bigl(u^{(n),s}_{j}-\widehat{m}^s_{j}\bigr)\otimes \bigl(u^{(n),s}_{j}-\widehat{m}^s_{j}\bigr) =\bigl(I - \msK(\widehat{\Sigma}_{j}^s ) H \bigr) \widehat{\Sigma}_{j}^s.
\end{split}
\end{align*}
\vspace{-0.4cm}
\STATE {\bf Output}: Means and covariances  $\{\widehat{m}_j^s, \widehat{C}_j^s\}_{j=1}^\infty.$
\end{algorithmic}
\end{algorithm}

We now state our second main result, which establishes long-time accuracy of ensemble Kalman filters with a surrogate model. We recall that $\eps$ represents the level of observation noise and $\delta$ the accuracy of the surrogate model. The proof can be found in Section \ref{sec:proofmain2}.

\begin{theorem}\label{th:main2}
   Suppose that Assumptions \ref{assumption:ball and squeezing} and \ref{assumption:surrogate model assumption} hold and that $N \ge 6k$. Then, if the inflation parameter $a>0$ is sufficiently large, there is a constant $\mathsf{C}^s$ independent of $\eps$ and $\delta$ such that
    \begin{align*}
        \limsup_{j\to \infty}\E \|\widehat{m}^s_j-u_j\|\leq \mathsf{C}^s(\eps+\delta).
    \end{align*}
\end{theorem}

\section{Filter accuracy via mean-field theory}\label{sec:proofmain1}
In this section, we prove the first main result, Theorem \ref{th:main1}. Our analysis will leverage new accuracy bounds for an idealized, non-implementable algorithm that uses population moments, rather than moments computed through an ensemble. Specifically, we consider the following Gaussian projected filter, which arises from a mean-field limit ($N\to \infty)$ of the ensemble Kalman filter in Algorithm \ref{algorithm:finite ensemble}. Notice that this idealized algorithm does not propagate an ensemble of particles, and covariance inflation is enforced by adding a positive definite operator $Q = aP$ to the prediction covariance.

\begin{algorithm}[H]
\caption{Gaussian Projected Filter \label{algorithm:mean field}}
\begin{algorithmic}[1]
\STATE {\bf Initialization}:  $m_0, C_0.$ 
\STATE For $j =  1, 2 \ldots$ do the following:
\STATE {\bf \index{prediction}Prediction}: 
\vspace{-0.3cm}
\begin{align*}
\begin{split}
\mu_{j} &= \mathbb{E}_{u\sim \Nc(m_{j-1},C_{j-1})}\bigl[\Psi(P_{\B_V} u)\bigr],  \\
\Sigma_{j} &= \mathbb{E}_{u\sim \Nc(m_{j-1},C_{j-1})}\bigl[\left(\Psi(P_{\B_V}u)-\mu_{j} \right)\otimes\left(\Psi(P_{\B_V}u)-\mu_{j} \right)\bigr].
\end{split}
\end{align*}
\vspace{-0.5cm}
\STATE{{{\bf \index{analysis}Analysis}}}: 
\vspace{-0.3cm}
	\begin{align*}
	\begin{split}
        m_{j}&=\mu_{j}+ \msK(\Sigma_{j} + Q)\left(y_{j} -H\mu_{j}\right),\\
        C_{j} &= \bigl( I - \msK(\Sigma_{j} + Q) H \bigr) (\Sigma_{j} + Q).
\end{split}
\end{align*}
\vspace{-0.4cm}
\STATE {\bf Output}: Means and covariances  $\{m_j, C_j\}_{j=1}^\infty.$
\end{algorithmic}
\end{algorithm}

 Our first lemma shows that the squeezing property implies that the trace of the analysis covariance operators generated by Algorithms \ref{algorithm:finite ensemble} and \ref{algorithm:mean field} are controlled at the noise level $\eps$ in the long-time asymptotic.
\begin{lemma}\label{lemma:covariance trace bound}
    Under Assumption \ref{assumption:ball and squeezing}, the covariances output by Algorithm \ref{algorithm:mean field} satisfy
    \begin{align}\label{eq:mean field covariance bound}
        \Tr \left(C_j \right)\leq \alpha^j \Tr \left(C_0 \right)+(1+\alpha \beta) \frac{1-\alpha^j}{1-\alpha}\eps^2\Tr(R),
    \end{align}
    and those output by Algorithm \ref{algorithm:finite ensemble} satisfy
    \begin{align}\label{eq:finite ensemble covariance bound}
        \Tr \bigl(\widehat{C}_j \bigr)\leq \alpha^j \Tr \bigl(\widehat{C}_0 \bigr)+(1+\alpha \beta) \frac{1-\alpha^j}{1-\alpha}\eps^2\Tr(R). 
        \end{align}
                \begin{proof}

Since $P$ is an orthogonal projection and $Q = aP,$ we have
            \begin{align}\label{eq:auxbound}
            \begin{split}
                    \Tr(C_{j})&=\Tr(P C_{j}P^*)+\Tr \bigl((I-P)C_{j}(I-P)^* \bigr)\\
            & \le \Tr(P C_{j}P^*)+\Tr \bigl((I-P)\Sigma_j (I-P)^* \bigr).
            \end{split}
            \end{align}
We next bound in turn the two terms in the right-hand side. 

For the first term,  let $\{\phi_i\}_{i=1}^k$ be an orthonormal basis for the range of $P$, which we can extend to a basis for $\H$ via Gram-Schmidt by the separability of $\H$. We then compute that
            \begin{align}
            \begin{split}\label{eq:firstterm}
\Tr\left(PC_{j}P^*\right)&=\sum_{i=1}^{\infty}\langle \phi_i, PC_{j}P^* \phi_i \rangle
                =\sum_{i=1}^{k}\langle \phi_i, PC_{j}P^* \phi_i \rangle\\
                &=\sum_{i=1}^k\langle H \phi_i, HC_jH^* H\phi_i \rangle_{\R^k}= \Tr\left(HC_jH^* \right) \le \eps^2\Tr\left(R\right),
            \end{split}
            \end{align}
           where we have used that $\{H\phi_i\}_{i=1
            }^k$ is an orthonormal basis for $\R^k$ and that $  \epsilon^2 R - HC_{j}H^*$ is positive semi-definite in the last inequality.           
                
         For the second term in the right-hand side of \eqref{eq:auxbound},   note that
            \begin{align}\label{eq:secondterm}
            \begin{split}
                \Tr \bigl( (I-P)\Sigma_j(I-P)^* \bigr)&=\mathbb{E}_{\Nc(m_{j-1},C_{j-1})} \Bigl[ \bigl\|(I-P)(\Psi(P_{\B_V} u)-\mu_j) \bigr\|^2 \Bigr]\\
                & \hspace{-1cm} \leq \mathbb{E}_{\Nc(m_{j-1},C_{j-1})}\Bigl[ \bigl\|(I-P)\bigl(\Psi(P_{\B_V} u)-\Psi(P_{\B_V}m_{j-1}) \bigr) \bigr\|^2 \Bigr]\\ 
                & \hspace{-1cm} \leq\mathbb{E}_{\Nc(m_{j-1},C_{j-1})} \Bigl[ \alpha \left( \bigl\|P_{\B_V} u-P_{\B_V} m_{j-1} \bigr\|^2+\beta \bigl\|P(P_{\B_V} u-P_{\B_V} m_{j-1})\bigr\|^2 \right) \Bigr]\\
                & \hspace{-1cm} \leq\mathbb{E}_{\Nc(m_{j-1},C_{j-1})} \Bigl[ \alpha \left(\|u-m_{j-1}\|^2+\beta \|P(u-m_{j-1})\|^2 \right) \Bigr]\\
                & \hspace{-1cm} =\alpha \Tr \left( C_{j-1}\right)+\alpha \beta \Tr \left(P C_{j-1}P^* \right).
                \end{split}
            \end{align}

        Combining \eqref{eq:auxbound}, \eqref{eq:firstterm}, and \eqref{eq:secondterm}, we deduce that 
         \begin{align*}
                \Tr(C_{j})& \leq \eps^2\Tr(R)+\alpha \Tr(C_{j-1})+\alpha\beta \Tr(PC_{j-1}P^*)\\
                & \leq (1+\alpha \beta)\eps^2\Tr (R)+\alpha \Tr (C_{j-1}).
            \end{align*}
      Applying the discrete Gronwall lemma \cite[Theorem 1.19]{sanz2023inverse} yields \eqref{eq:mean field covariance bound}. 
        
            To prove  \eqref{eq:finite ensemble covariance bound}, we proceed similarly. Since $Q=aP$, it holds that $(I-P)\xi_{j}^{(n)}=0$ almost surely, from which we have 
            \begin{align*}
                \Tr \bigl( (I-P)\widehat{\Sigma}_j(I-P)^*\bigr)&=\frac{1}{N-1}\sum_{n=1}^N \Bigl\|(I-P)(v_{j}^{(n)}-\widehat{\mu}_j) \Bigr\|^2\\
                & \hspace{-1cm} =\frac{1}{N-1}\sum_{n=1}^N \Bigl\|(I-P)\bigl(\Psi(P_{\B_V}u^{(n)}_{j-1})-\frac{1}{N}\sum_{n=1}^N\Psi(P_{\B_V} u^{(n)}_{j-1})\bigr) \Bigr\|^2\\
                & \hspace{-1cm}\leq \frac{1}{N-1}\sum_{n=1}^N \Bigl\|(I-P)\bigl(\Psi(P_{\B_V}u^{(n)}_{j-1})-\Psi(P_{\B_V}\widehat{m}_{j-1})\bigr) \Bigr\|^2\\
                & \hspace{-1cm} \leq \frac{1}{N-1}\sum_{n=1}^N\alpha \left( \bigl\|P_{\B_V}u^{(n)}_{j-1}-P_{\B_V}\widehat{m}_{j-1} \bigr\|^2+\beta \bigl\|P(P_{\B_V}u^{(n)}_{j-1}-P_{\B_V}\widehat{m}_{j-1}) \bigr\|^2 \right)\\
                & \hspace{-1cm} \leq \frac{1}{N-1}\sum_{n=1}^N\alpha \left( \bigl\|u^{(n)}_{j-1}-\widehat{m}_{j-1} \bigr\|^2+\beta \bigl\|P(u^{(n)}_{j-1}-\widehat{m}_{j-1}) \bigr\|^2 \right)\\
                & \hspace{-1cm} =\alpha \Tr\bigl( \widehat{C}_{j-1}\bigr)+ \alpha \beta \Tr \bigl( P \widehat{C}_{j-1}P^*\bigr).
            \end{align*}
            The remainder of the proof proceeds exactly as for the mean-field version.
        \end{proof}
\end{lemma}
The next theorem establishes long-time accuracy of the idealized mean-field Algorithm \ref{algorithm:mean field}.
\begin{theorem}\label{theorem:general mean field filter accuracy}
Suppose that Assumption \ref{assumption:ball and squeezing} holds. Then, if the inflation parameter $a>0$ is sufficiently large, there is a constant $\mathsf{C}_1$ independent of $\eps$ such that
\begin{align*}
    \limsup_{j\to \infty} \mathbb{E}\|m_j-u_j\|\leq \mathsf{C}_1\eps.
\end{align*}
    \begin{proof}
        From Lemma \ref{lemma:covariance trace bound}, there exists a time $j_*$ such that for all $j\geq j_*$, $\Tr(C_j)\leq \left(1+\frac{1+\alpha \beta}{1-\alpha}\Tr(R) \right)\eps^2$. Throughout the rest of the proof, we assume that $j \ge j_*.$ We assume without loss of generality that $u_0\in \B$, since if $u_0\notin \B$, there exists a time $j$ such that $u_j\in \B$, and we take $j_*>j.$
        Next, we write 
        \begin{align}\label{eq:signal}
        \begin{split}
            u_j &= \Psi(u_{j-1}) \\
            & = \bigl(I-\msK(\Sigma_j+Q)H \bigr)\Psi(u_{j-1})+\msK(\Sigma_j+Q)H \Psi(u_{j-1}),  \\
        \end{split}    
        \end{align}
        and introduce the auxiliary variable
        \begin{align}\label{eq:auxvariablew}
        \begin{split}
            w_j &:=\bigl(I-\msK(\Sigma_j+Q)H \bigr)\Psi(P_{\B_V}w_{j-1})+\msK(\Sigma_j+Q)y_j, \\
            & = \bigl(I-\msK(\Sigma_j+Q)H \bigr)\Psi(P_{\B_V}w_{j-1})+\msK(\Sigma_j+Q)H\Psi(u_{j-1}) + \eps \msK(\Sigma_j+Q)H\eta_j.
         \end{split}
        \end{align}
        Subtracting \eqref{eq:signal} from \eqref{eq:auxvariablew}, we obtain that
        \begin{align*}
            w_j-u_j&=\bigl(I- \msK(\Sigma_j+Q)H\bigr)\bigl( \Psi(P_{B_V}w_{j-1})-\Psi(u_{j-1})\bigr)+\eps \msK(\Sigma_j+Q)\eta_j\\
            &=\bigl(I-\msK(\Sigma_j+Q)H \bigr)P\bigl(\Psi(P_{B_V}w_{j-1}) -\Psi(u_{j-1})\bigr)
            \\ &+\bigl(I-P \bigr)\bigl(\Psi(P_{B_V}w_{j-1}) -\Psi(u_{j-1})\bigr)+ \eps \msK(\Sigma_j+Q)\eta_j,
        \end{align*}
        where we used that $H(I-P) = 0.$ Taking the $V$ norm of both sides and using that $P^2 = P$ yields
        \begin{align*}
            V\bigl(w_j-u_j \bigr)&\leq \bigl\|P-\msK(\Sigma_j+Q)H \bigr\|_{V,V}V\Bigl(P\bigl(\Psi(P_{B_V}w_{j-1}) -\Psi(u_{j-1})\bigr) \Bigr)\\
            &  +V\Bigl(\bigl(I-P \bigr)\bigl(\Psi(P_{B_V}w_{j-1}) -\Psi(u_{j-1})\bigr) \Bigr)+\eps \bigl\|\msK(\Sigma_j+Q)H \bigr\|_{V,V}V(H^*\eta_j)\\
            &\leq \bigl\|P-\msK(\Sigma_j+Q)H \bigr\|_{V,V}L^{1/2}V\bigl(w_{j-1} -u_{j-1}\bigr)\\
            & +\alpha^{1/2}V\bigl(w_{j-1}-u_{j-1}\bigr)+\eps \bigl\|\msK(\Sigma_j+Q)H \bigr\|_{V,V}V(H^*\eta_j),
        \end{align*}
        where the last inequality follows by Assumption \ref{assumption:ball and squeezing}. Now, note that
        \begin{align*}
            \bigl\|P-\msK(\Sigma_j+Q)H \bigr\|_{V,V}&=\sup_{V(v)=1}V\Bigl( \bigl(P-\msK(\Sigma_j+Q)H \bigr)Pv\Bigr)\\
         & \hspace{-2.4cm} \leq \sup_{V(v)=1}\sqrt{2}\Bigl( \bigl\| \bigl(P-\msK(\Sigma_j+Q)H \bigr)Pv \bigr\|+ \bigl\|(I-P) \bigl(P-\msK(\Sigma_j+Q)H \bigr)Pv \bigr\| \Bigr)\\
        &  \hspace{-2.4cm} \leq \sup_{V(v)=1}\sqrt{2}\Bigl( \bigl\|(P-\msK(\Sigma_j+Q)H)Pv \bigr\|+\sup_{V(v)=1} \bigl\|(I-P)\msK(\Sigma_j+Q)Hv \bigr\|\Bigr).
        \end{align*}
        The first of the terms on the right-hand side can be bounded as
        \begin{align*}
            \bigl\| \bigl(P-\msK(\Sigma_j+Q)H \bigr)Pv \bigr\|&\leq \Bigl\|I-H(\Sigma_j+Q)H^*\bigl(H(\Sigma_j+Q)H^*+\eps^2R \bigr)^{-1} \Bigr\|_{op}\|Pv\|\\
            &\leq \eps^2\frac{\|R\|_{op}}{\lambda_{\min }(H\Sigma_jH^*+aI)}\|Pv\|\\
            & \leq \eps^2\frac{\|R\|_{op}}{a}\|Pv\|
             \leq \eps^2\frac{\|R\|_{op}}{a}V(v).
        \end{align*}
        For the second term, we have
        \begin{align*}
           \bigl\|(I-P)\msK(\Sigma_j+Q)Hv \bigr\| & \leq \bigl\|(I-P)(\Sigma_j+Q) \bigr\|_{op} \Bigl\|\bigl(H(\Sigma_j+Q)H^*+\eps^2R \bigr)^{-1} \Bigr\|_{op}\|Hv\|\\
            & \leq \frac{\|\Sigma_j\|_{op}}{\lambda_{\min}\bigl(H(\Sigma_j+Q)H^* \bigr)}\|Hv\|\\
            & \leq \frac{L\|C_{j-1}\|_{op}}{a}\|Hv\|
             \leq \frac{L\Tr(C_{j-1})}{a}\|Hv\|\\
            & \leq \frac{L\left(1+\frac{1+\alpha \beta}{1-\alpha}\Tr(R) \right) \eps^2}{a}\|Hv\|
             \leq \frac{L\left(1+\frac{1+\alpha \beta}{1-\alpha}\Tr(R) \right)\eps^2}{a}V(v).
        \end{align*}
        Consequently, 
        \begin{align}\label{eq:mf P minus Kalman gain norm bound}
            \bigl\|P-\msK(\Sigma_j+Q)H \bigr\|_{V,V}\leq \frac{\eps^2}{a}\biggl(\|R\|_{op}+L\Bigl(1+\frac{1+\alpha \beta}{1-\alpha}\Tr(R) \Bigr) \biggr).
        \end{align}
       This implies that
        \begin{align}\label{eq:kalman gain bound}
            \|\msK(\Sigma_j+Q)H\|_{V,V}\leq 1+\frac{\eps^2}{a}\biggl(\|R\|_{op}+L\Bigl(1+\frac{1+\alpha \beta}{1-\alpha}\Tr(R) \Bigr) \biggr).
        \end{align}
        By taking $a$ large enough to guarantee that $ \bigl\|P-\msK(\Sigma_j+Q)H \bigr\|_{V,V}L^{1/2}\leq \frac{1-\alpha^{1/2}}{2}$ we obtain 
        \begin{align*}
             V \bigl(w_j-u_j \bigr)&\leq \alpha_*V \bigl(w_{j-1}-u_{j-1}\bigr) \\
             &+\eps \Biggl(1+\frac{\eps^2}{a}\biggl(\|R\|_{op}+L\Bigl(1+\frac{1+\alpha \beta}{1-\alpha}\Tr(R) \Bigr) \biggr) \Biggr)V(H^*\eta_j),
        \end{align*}
        where $\alpha_*:=\frac{1+\alpha^{1/2}}{2}<1.$ Taking the expectation of both sides of the inequality we have
        \begin{align*}
            \E   V \bigl(w_j-u_j \bigr) &\leq \alpha_*\E  V \bigl(u_{j-1}-w_{j-1} \bigr)  \\ & +\eps \Biggl(1+\frac{\eps^2}{a}\biggl(\|R\|_{op}+L\Bigl(1+\frac{1+\alpha \beta}{1-\alpha}\Tr(R) \Bigr) \biggr)+ \sqrt{2\Tr (R)}\Biggr).
        \end{align*}
        Applying the discrete Gronwall inequality yields
        \begin{align}\label{eq: gronwall bound 1}
            \E  V \bigl(w_j-u_j \bigr) \leq \frac{c_1\eps}{1-\alpha_*}\left(1-\alpha_*^{j-j_*} \right)+\alpha_*^{j-j_*}\E  V\bigl(u_{j_*}-w_{j_*} \bigr),
        \end{align}
        where $c_1=\left(1+\frac{\eps^2}{a}\left(\|R\|_{op}+L\left(1+\frac{1+\alpha \beta}{1-\alpha}\Tr(R) \right) \right)+ \sqrt{2\Tr (R)}\right).$
        
        We now consider the difference
        \begin{align*}
            m_j-w_j&=\bigl(I-\msK(\Sigma_j+Q)H \bigr)\bigl(\mu_j-\Psi(P_{\B_V}w_{j-1}) \bigr)\\
            &=\bigl(I-\msK(\Sigma_j+Q)H \bigr)\bigl(\mu_j-\Psi(P_{\B_V}m_{j-1}) \bigr)  \\
            &+\bigl(I-\msK(\Sigma_j+Q)H \bigr)\bigl(\Psi(P_{\B_V}m_{j-1})-\Psi(P_{\B_V}w_{j-1} ) \bigr).
        \end{align*}
        Carrying out the same argument as above, we have that
        \begin{align*}
            V \bigl(m_j-w_j \bigr)&\leq V\Bigl( \bigl(I-\msK(\Sigma_j+Q)H \bigr)\bigl(\mu_j-\Psi(P_{\B_V}m_{j-1}) \bigr) \Bigr)+ \alpha_* V\bigl(m_{j-1}-w_{j-1} \bigr)\\
            & \hspace{-1cm} \leq \|\left(I-\msK(\Sigma_j+Q)H \right)\|_{V,V}V\bigl(\mu_j-\Psi(P_{\B_V}m_{j-1})\bigr)+ \alpha_* V\bigl(m_{j-1}-w_{j-1} \bigr)\\
            & \hspace{-1cm}  \leq \Biggl(1+\frac{\eps^2}{a}\Biggl(\|R\|_{op}+L\Bigl(1+\frac{1+\alpha \beta}{1-\alpha}\Tr(R) \Bigr) \biggr) \Biggr)V\bigl(\mu_j-\Psi(P_{\B_V}m_{j-1})\bigr)  \\
            & \hspace{-1cm}  + \alpha_* V\bigl(m_{j-1}-w_{j-1} \bigr).
        \end{align*}
        We now bound $V\bigl(\mu_j-\Psi(P_{\B_V}m_{j-1})\bigr):$
        \begin{align*}
            V\bigl(\mu_j-\Psi(P_{\B_V}m_{j-1})\bigr) &=V\left(\mathbb{E}_{u\sim \Nc(m_{j-1},C_{j-1})}\bigl[\Psi(P_{\B_V}u)-\Psi(P_{\B_V}m_{j-1})\bigr] \right)\\
            &\leq \Bigl( \E_{u\sim \Nc(0,C_{j-1})}\Bigl[V\Bigl(\Psi\bigl(P_{\B_V}(u-m_{j-1})\bigr)-\Psi \bigl(P_{\B_V}m_{j-1}\bigr) \Bigr)\Bigr]\Bigr)^{1/2}\\
            &\leq L^{1/2}\left( \E_{u\sim \Nc(0,C_{j-1})}\Bigl[V\Bigl(P_{\B_V}(u-m_{j-1})-P_{\B_V}m_{j-1}\right)\Bigr]\Bigr)^{1/2}\\
            &\leq L^{1/2}\left( \E_{u\sim \Nc(0,C_{j-1})}\bigl[V\left(u-m_{j-1}\right)\bigr]\right)^{1/2}\\
            & \leq L^{1/2}\bigl(\Tr(C_{j-1})+\Tr(HC_{j-1}H^*) \bigr)^{1/2}\\
            & \leq \eps L^{1/2}\Biggl(1+\biggl(\frac{1+\alpha \beta}{1-\alpha}+1\biggr)\Tr(R) \Biggr)^{1/2}.
        \end{align*}
        Consequently, 
        \begin{align*}
            \E V\bigl(m_j-w_j \bigr)& \leq \alpha_*\E V\bigl(m_{j-1}-w_{j-1} \bigr)+c_2\eps,
        \end{align*}
        where $c_2=L^{1/2}\biggl(1+\Bigl(\frac{1+\alpha \beta}{1-\alpha}+1\Bigr)\Tr(R) \biggr)^{1/2}\Biggl(1+\frac{\eps^2}{a}\biggl(\|R\|_{op}+L\Bigl(1+\frac{1+\alpha \beta}{1-\alpha}\Tr(R) \Bigr) \biggr) \Biggr).$
        Therefore, the discrete Gronwall inequality gives
        \begin{align}\label{eq: gronwall bound 2}
            \E  V \bigl(m_j-w_j \bigr)  \leq \frac{c_2\eps}{1-\alpha_*}\left(1-\alpha_*^{j-j_*} \right)+\alpha_*^{j-j_*}\E V \bigl(m_{j_*}-w_{j_*} \bigr).
        \end{align}
        Combining \eqref{eq: gronwall bound 1} and \eqref{eq: gronwall bound 2}, we have
        \begin{align*}
            \E  V\bigl(m_j-u_j \bigr)  &\leq \frac{(c_1+c_2)\eps}{1-\alpha_*}\left(1-\alpha_*^{j-j_*} \right)  +\alpha_*^{j-j_*}\E V \bigl(u_{j_*}-w_{j_*}\bigr) +\alpha_*^{j-j_*}\E  V \bigl(m_{j_*}-w_{j_*} \bigr).
        \end{align*}
       Since $\alpha_*\in (0,1)$ it follows that, for $c_3=\frac{c_1+c_2}{1-\alpha^*},$
       \begin{align*}
           \limsup_{j\to \infty}\E V\bigl(m_j-u_j \bigr)\leq c_3\eps.
       \end{align*}
       The equivalence of the $V(\cdot)$ and $\|\cdot \|$ norms concludes the proof.
    \end{proof}
\end{theorem}

We now show that the means output by the square-root filter Algorithm \ref{algorithm:finite ensemble} and those output by the idealized mean-field Algorithm \ref{algorithm:mean field} remain close in the long-time asymptotic. The proof of our first main result, Theorem \ref{th:main1}, follows directly by combining Theorem \ref{theorem:general mean field filter accuracy} 
and the next result: 
\begin{theorem}\label{theorem:general finite ensemble filter accuracy}
Suppose that Assumption \ref{assumption:ball and squeezing} holds and that $N\geq 6k$. Then, if the inflation parameter $a>0$ is sufficiently large, there is a constant $\mathsf{C}_2$ independent of $\eps$ such that
\begin{align}
    \limsup_{j\to \infty} \E \|\widehat{m}_j-m_j\|\leq \mathsf{C}_2 \eps.
\end{align}
\begin{proof}
    We assume that $a$ is chosen to be at least as large as in Theorem \ref{theorem:general mean field filter accuracy}. From Theorem \ref{theorem:general mean field filter accuracy} and Lemma \ref{lemma:covariance trace bound}, there exists a time $j_*$ such that for all $j>j_*$, $\Tr(\widehat{C}_j)\leq \left(1+\frac{1+\alpha \beta}{1-\alpha}\Tr(R) \right)\eps^2$, $\Tr(C_j)\leq \left(1+\frac{1+\alpha \beta}{1-\alpha}\Tr(R) \right)\eps^2$, and $\E V(m_j-u_j)\leq \max\{1,2\mathsf{C}_1\eps\}.$ For $j>j_*,$ we consider the difference
    \begin{align*}
        \widehat{m}_j-m_j&=\bigl(I-\mathscr{K}(\widehat{\Sigma}_j)H \bigr)\widehat{\mu}_j-\bigl(I-\mathscr{K}(\Sigma_j+Q)H \bigr)\mu_j+\bigl(\mathscr{K}(\widehat{\Sigma}_j)-\mathscr{K}(\Sigma_j)\bigr)y_j \\
        &\hspace{-1cm}  =\bigl(I- \mathscr{K}(\Sigma_j+Q)H\bigr)\bigl(\mu_j-\widehat{\mu_j} \bigr) 
        +\bigl(I-\mathscr{K}(\Sigma_j+Q)H-I+\mathscr{K}(\widehat{\Sigma}_j)H \bigr)\bigl(P\widehat{\mu}_j-H^*y_j \bigr).
    \end{align*}
    Taking the $V$ norm of both sides, we have
    \begin{align}\label{eq:finite sample inequality 1}
    \begin{split}
        V\left(\widehat{m}_j-m_j \right) &\leq V\Bigl( \bigl(I- \mathscr{K}(\Sigma_j+Q)H\bigr)\bigl(\mu_j-\widehat{\mu_j} \bigr)\Bigr) \\
        &+V\Bigl(\bigl(\mathscr{K}(\widehat{\Sigma}_j)H -\mathscr{K}(\Sigma_j+Q)H\bigr)\bigl(P\widehat{\mu}_j-H^*y_j \bigr) \Bigr).
        \end{split}
    \end{align}
    We begin by bounding the first term on the right-hand side of \eqref{eq:finite sample inequality 1}
    \begin{align*}
        V\Bigl( \bigl(I- \mathscr{K}(\Sigma_j+Q)H\bigr)\bigl(\mu_j-\widehat{\mu_j} \bigr)\Bigr)
        &\leq V\Bigl( \bigl(I- \mathscr{K}(\Sigma_j+Q)H\bigr)\bigl(\mu_j-\Psi(P_{B_V}m_{j-1}) \bigr)\Bigr)\\
        & \hspace{-2.5cm} +V\Bigl( \bigl(I- \mathscr{K}(\Sigma_j+Q)H\bigr)\bigl(\Psi(P_{B_V}m_{j-1})-\Psi(P_{B_V}\widehat{m}_{j-1}) \bigr)\Bigr)\\
        & \hspace{-2.5cm} +V\Biggl( \bigl(I- \mathscr{K}(\Sigma_j+Q)H\bigr)\biggl(\Psi(P_{B_V}\widehat{m}_{j-1})-\frac{1}{N}\sum_{n=1}^N\Psi \bigl(P_{B_V}u_{j-1}^{(n)} \bigr) \biggr)\Biggr)\\
        & \hspace{-2.5cm} +V\left( \bigl(I- \mathscr{K}(\Sigma_j+Q)H\bigr)\frac{1}{N}\sum_{n=1}^N\xi^{(n)}_j\right).
    \end{align*}
  Recall that $\|I-\msK(\Sigma_j+Q)H\|_{V,V}\leq c_1:=1+\frac{\eps^2}{a}\biggl(\|R\|_{op}+L\Bigl(1+\frac{1+\alpha \beta}{1-\alpha}\Tr(R) \Bigr) \biggr).$ By the same arguments as in the proof of Theorem \ref{theorem:general mean field filter accuracy}, we have 
    \begin{align*}
    V\Bigl( \bigl(I- \mathscr{K}(\Sigma_j+Q)H\bigr)\bigl(\mu_j-\Psi(P_{\B_V}m_{j-1}) \bigr)\Bigr)\leq c_2\eps ,
    \end{align*}
    and 
    \begin{align*}
        &V\Bigl( \bigl(I- \mathscr{K}(\Sigma_j+Q)H\bigr)\bigl(\Psi(P_{\B_V}m_{j-1})-\Psi(P_{\B_V}\widehat{m}_{j-1}) \bigr)\Bigr) \\
        & \hspace{5cm} \leq \left( \frac{\eps^2}{a}c_3+\alpha^{1/2} \right)V\left(m_{j-1}-\widehat{m}_{j-1} \right),
    \end{align*}
    with $c_2=c_1L^{1/2}\biggl(1+\Bigl(1+\frac{1+\alpha \beta}{a-\alpha}\Bigr)\Tr(R) \biggr)^{1/2}$ and $c_3=\|R\|_{op}+L\left(1+\frac{1+\alpha \beta}{1-\alpha}\Tr(R)\right).$ Moreover,
    \begin{align*}
        &V\Biggl( \bigl(I- \mathscr{K}(\Sigma_j+Q)H\bigr)\Bigl(\Psi(P_{\B_V}\widehat{m}_{j-1})-\frac{1}{N}\sum_{n=1}^N \Psi(P_{\B_V}u_{j-1}^{(n)}) \Bigr)\Biggr)\\
        & \leq c_1V\left( \Psi(P_{\B_V}\widehat{m}_{j-1})-\frac{1}{N}\sum_{n=1}^N \Psi(P_{\B_V}u_{j-1}^{(n)})\right)\\
        & \leq c_1 \left(\frac{L}{N}\sum_{n=1}^NV\left( u_{j-1}^{(n)}-\widehat{m}_{j-1}\right) \right)^{1/2} \leq  c_1L^{1/2}\left(\Tr(\widehat{C}_{j-1})+\Tr(H\widehat{C}_{j-1}H^*)\right)^{1/2} \\
        & \leq\eps c_1L^{1/2}\left(1+\left(\frac{1+\alpha \beta}{1-\alpha}+1\right)\Tr(R) \right)^{1/2}
        \leq c_4 \eps,
    \end{align*}
    where $c_4=c_1L^{1/2}\left(1+\left(\frac{1+\alpha \beta}{1-\alpha}+1\right)\Tr(R) \right)^{1/2}$.
    Finally, we have
    \begin{align*}
        \E V\left( \bigl(I- \mathscr{K}(\Sigma_j+Q)H\bigr)\frac{1}{N}\sum_{n=1}^N\xi^{(n)}_j\right)&\leq 2 \bigl\|P-\mathscr{K}(\Sigma_j+Q)H \bigr\|_{V,V} \, \E \Bigl\|\frac{1}{N}\sum_{n=1}H\xi_{j}^{(n)} \Bigr\|\\
        & \hspace{-2cm}  \leq  2\frac{\eps^2}{a}\biggl(\|R\|_{op}+L\Bigl(1+\frac{1+\alpha \beta}{1-\alpha}\Tr(R) \Bigr) \biggr)\E \Bigl\|\frac{1}{N}\sum_{n=1}H\xi_{j}^{(n)} \Bigr\|\\ 
        & \hspace{-2cm}  \overset{\text{(i)}}{\leq} 2\frac{\eps^2}{a}\biggl(\|R\|_{op}+L\Bigl(1+\frac{1+\alpha \beta}{1-\alpha}\Tr(R) \Bigr) \biggr)\sqrt{\frac{\Tr(aI)}{N}}\\
        & \hspace{-2cm}  \overset{\text{(ii)}}{\leq }\eps^2a^{-1/2}\biggl(\|R\|_{op}+L\Bigl(1+\frac{1+\alpha \beta}{1-\alpha}\Tr(R) \Bigr) \biggr),
    \end{align*}
    where (i) uses Lemma B.9 in \cite{al2024ensemble} and (ii) follows from the assumption that $N\geq 6k$. Consequently,
    \begin{align}\label{eq:finite sample bound 2}
        \E V\Bigl(\bigl(I- \mathscr{K}(\Sigma_j+Q)H\bigr)\bigl(\mu_j-\widehat{\mu_j} \bigr) \Bigr)\leq \left(\frac{\eps^2}{a}c_3+\alpha^{1/2} \right) \E V\bigl(\widehat{m}_{j-1}-m_{j-1} \bigr)+c_5\eps,
    \end{align}
    where
    \begin{align*}
        c_5&:=c_2+c_4+\eps a^{-1/2}\biggl(\|R\|_{op}+L\Bigl(1+\frac{1+\alpha \beta}{1-\alpha}\Tr(R) \Bigr) \biggr).
    \end{align*}
    We now bound the second term in \eqref{eq:finite sample inequality 1}:
    \begin{align*}
        &V\Bigl(\bigl(\mathscr{K}(\widehat{\Sigma}_j)H -\mathscr{K}(\Sigma_j+Q)H\bigr)\bigl(P\widehat{\mu}_j-H^*y_j \bigr) \Bigr)\\
        &\leq V\Bigl(\bigl(P-\mathscr{K}(\widehat{\Sigma}_j)H\bigr)\bigl(P\widehat{\mu}_j-H^*y_j \bigr) \Bigr)+V\Bigl(\bigl(P-\mathscr{K}(\Sigma_j+Q)H\bigr)\bigl(P\widehat{\mu}_j-H^*y_j \bigr) \Bigr)\\
        & \leq \left(\|P-\mathscr{K}(\Sigma_j+Q)H\|_{V,V}V+\|P-\mathscr{K}(\widehat{\Sigma}_j)H\|_{V,V}\right)V\bigl(P\widehat{\mu}_j-H^*y_j \bigr)\\
        & \leq \frac{\eps^2}{\lambda_{\min}(H\Sigma_jH^*+aI)}c_3V\bigl(P\widehat{\mu}_j-H^*y_j \bigr)
         +  \frac{\eps^2}{\lambda_{\min}(H\widehat{\Sigma}_jH^*)}c_3V\bigl(P\widehat{\mu}_j-H^*y_j \bigr).
    \end{align*}
    Taking expectations,
    \begin{align*}
        &\E  V\Bigl(\bigl(\mathscr{K}(\widehat{\Sigma}_j)H -\mathscr{K}(\Sigma_j+Q)H\bigr)\bigl(P\widehat{\mu}_j-H^*y_j \bigr) \Bigr)\\
        & \leq \frac{\eps^2}{a}c_3\E V\bigl(P\widehat{\mu}_j-H^*y_j \bigr)+ \eps^2c_3\E \biggl[\frac{1}{\lambda_{\min}(H\widehat{\Sigma}_jH^*)} V\left(P\widehat{\mu}_j-H^*y_j \right)\biggr]\\
        & \leq \frac{\eps^2}{a}c_3\E V\bigl(P\widehat{\mu}_j-H^*y_j \bigr)+ \eps^2c_3\E \Biggl[\E \biggl[\frac{1}{\lambda_{\min}(H\widehat{\Sigma}_jH)}V\bigl(P\widehat{\mu}_j-H^*y_j\bigr) \Big\vert \{u_{j-1}^{(n)}\}_{n=1}^N\biggr]\Biggr].
    \end{align*}
    Denoting $\mathbb{E}^\mathsf{U_{j-1}} [\, \cdot \,] := \mathbb{E} \Bigl[ \, \cdot \, \bigl| \{u_{j-1}^{(n)}\}_{n=1}^N \Bigr]$, we then compute
    \begin{align*}
        &\E^\mathsf{U_{j-1}} \Biggl[\frac{1}{\lambda_{\min}(H\widehat{\Sigma}_jH)}V\left(P\widehat{\mu}_j-H^*y_j\right)\Biggr]\\ & 
        \leq \E^\mathsf{U_{j-1}} \Biggl[\frac{1}{\lambda_{\min}(H\widehat{\Sigma}_jH)}\Biggl(V\left(\frac{1}{N}\sum_{n=1}^NP\Psi(P_{\B_V}u_{j-1}^{(n)})-P\Psi(P_{\B_V}\widehat{m}_{j-1})\right)\\
        &+V\Bigl(P\Psi(P_{\B_V}\widehat{m}_{j-1})-P\Psi(P_{\B_V}m_{j-1})\Bigr)+V\Bigl(P\Psi(P_{\B_V}m_{j-1})-P\Psi(u_{j-1})\Bigr)\\
        &+\eps V\left(H^*\eta_j\right)+V\left(\frac{1}{N}\sum_{n=1}^N\xi_{j}^{(n)}\right)\Biggr)\Biggr]\\
        & \overset{\text{(i)}}{\leq} \E^\mathsf{U_{j-1}} \Biggl[\frac{1}{\lambda_{\min}(H\widehat{\Sigma}_jH)} \Biggr]\E^\mathsf{U_{j-1}} \Biggl[V\left(\frac{1}{N}\sum_{n=1}^N\Psi(P_{\B_V}u_{j-1}^{(n)})-\Psi(P_{\B_V}\widehat{m}_{j-1})\right)\\
        &+V\Bigl(\Psi(P_{\B_V}\widehat{m}_{j-1})-\Psi(P_{\B_V}m_{j-1})\Bigr)+V\Bigl(\Psi(P_{\B_V}m_{j-1})-\Psi(u_{j-1})\Bigr)+\eps V\left(H^*\eta_j\right) \Biggr]\\
        & + \E^\mathsf{U_{j-1}} \Biggl[\frac{1}{\lambda_{\min}(H\widehat{\Sigma}_jH)}V\left(\frac{1}{N}\sum_{n=1}^N\xi_j^{(n)}\right)\Biggr] \\
        &  \overset{\text{(ii)}}{\leq} \E^\mathsf{U_{j-1}} \Biggl[\frac{1}{\lambda_{\min}(H\widehat{\Sigma}_jH)}\Biggr]\E^\mathsf{U_{j-1}} \Biggl[V\left(\frac{1}{N}\sum_{n=1}^N\Psi(P_{\B_V}u_{j-1}^{(n)})-\Psi(P_{\B_V}\widehat{m}_{j-1})\right)\\
        &+V\Bigl(\Psi(P_{\B_V}\widehat{m}_{j-1})-\Psi(P_{\B_V}m_{j-1})\Bigr)+V\Bigl(\Psi(P_{\B_V}m_{j-1})-\Psi(u_{j-1})\Bigr)+\eps V\left(H^*\eta_j\right) \Biggr]\\
        & + \left(\E^\mathsf{U_{j-1}} \Biggl[\frac{1}{\lambda_{\min}(H\widehat{\Sigma}_jH)^2}\Biggr]  \E^\mathsf{U_{j-1}} \Biggl[ V^2\left(\frac{1}{N}\sum_{n=1}^N\xi_j^{(n)}\right)\Biggr]\right)^{1/2},
    \end{align*}
    
    where (i) uses conditional independence and (ii) uses Cauchy-Schwarz. If $a$ is chosen to guarantee that $\frac{10NL\left(1+\frac{1+\alpha \beta}{1-\alpha}\Tr(R) \right)\eps^2}{ak}\leq 1$, Lemma \ref{lemma:expected inverse eigenvalue} yields that
    \begin{align*}
        \E^\mathsf{U_{j-1}} \Biggl[\frac{1}{\lambda_{\min}(H\widehat{\Sigma}_jH)}\Biggr]\leq \frac{2C'}{a} \quad  \text{ and } \quad 
        \E^\mathsf{U_{j-1}} \Biggl[\frac{1}{\lambda_{\min}(H\widehat{\Sigma}_jH)^2} \Biggr]^{1/2}\leq \frac{\sqrt{2}C'}{a}.
    \end{align*}
   Note that both of these upper bounds are independent of $\{u_{j-1}^{(n)}\}_{n=1}^N.$ Consequently,
    \begin{align*}
       & \E \Biggl[ \E^\mathsf{U_{j-1}} \Bigl[\frac{1}{\lambda_{\min}(H\widehat{\Sigma}_jH)}V\left(P\widehat{\mu}_j-H^*y_j\right) \Bigr] \Biggr] \\ 
       &\leq \frac{2C'}{a}\E V\left(\frac{1}{N}\sum_{n=1}^N\Psi(P_{\B_V}u_{j-1}^{(n)})-\Psi(P_{\B_V}\widehat{m}_{j-1})\right)\\
    &+\E V\Bigl(\Psi(P_{\B_V}\widehat{m}_{j-1})-\Psi(P_{\B_V}m_{j-1})\Bigr) \\ &+\E V\Bigl(\Psi(P_{\B_V}m_{j-1})-\Psi(u_{j-1})\Bigr) 
        +\eps \E V\left(H^*\eta_j\right) 
         +  \E \biggl[V^2\Bigl(\frac{1}{N}\sum_{n=1}^N\xi_j^{(n)}\Bigr)\biggr]^{1/2}\\
        & \leq \frac{2C'}{a} \Biggl\{ L^{1/2}\left(\Tr\bigl(\widehat{C}_{j-1}\bigr)+\Tr\bigl(H\widehat{C}_{j-1}H^*\bigr)\right)^{1/2}\\
        & + L^{1/2}\E V\left(\widehat{m}_{j-1}-m_{j-1}\right)+L^{1/2}\E V\left(m_{j-1}-u_{j-1}\right)
         +\eps \sqrt{2 \Tr \left(R \right)}+\sqrt{\frac{\Tr(a I)}{N}}\Biggr\}.
    \end{align*}
    Since $j>j_*$, we have that 
    \begin{align}\label{eq:finite sample bound 3}
        \E V\bigl(\Psi(P_{\B_V}\widehat{m}_{j-1})-\Psi(P_{\B_V}m_{j-1})\bigr)\leq \frac{1}{a}\left(c_6+c_7\E V(\widehat{m}_{j-1}-m_{j-1})\right)+\frac{c_8}{a^{1/2}},
        \end{align}
        with $c_6=\eps 2\sqrt{2}C'\left(1+\frac{1+\alpha \beta}{1-\alpha}\Tr(R) \right)^{1/2}+2C'L^{1/2}\max\{2\mathsf{C}_1\eps,1\}+2C'\eps\sqrt{2\Tr(R)}$, $c_7=2C'L^{1/2}$, and $c_8=C'\sqrt{\frac{1}{3}}.$ By a similar computation, we have
        \begin{align} \label{eq:finite sample bound 4}
            \E V(P\widehat{\mu}_j-H^*y_j)\leq c_9+L^{1/2}\E V(\widehat{m}_{j-1}-m_{j-1})+c_{10}a^{1/2},
        \end{align}
        with $c_9=\eps \sqrt{2}\left(1+\frac{1+\alpha \beta}{1-\alpha}\Tr(R) \right)^{1/2}+L^{1/2}\max\{2\mathsf{C}_1\eps,1\}+\eps\sqrt{2\Tr(R)}$ and  $c_{10}=\sqrt{\frac{1}{6}}.$ Combining \eqref{eq:finite sample bound 3} and \eqref{eq:finite sample bound 4}, 
        \begin{align}\label{finite sample bound 5}
        \begin{split}
            &\E  V\Bigl(\bigl(\mathscr{K}(\widehat{\Sigma}_j)H -\mathscr{K}(\Sigma_j+Q)H\bigr)\bigl(P\widehat{\mu}_j-H^*y_j \bigr) \Bigr)\\ &\leq \frac{\eps^2}{a}c_3\left(c_9+L^{1/2}\E V(\widehat{m}_{j-1}-m_{j-1})+c_{10}a^{1/2}+c_6+c_7\E V(\widehat{m}_{j-1}-m_{j-1})+c_8a^{1/2}
       \right).
        \end{split}
        \end{align}
        Together, \eqref{eq:finite sample bound 2} and \eqref{eq:finite sample bound 4} yield
        \begin{align}
        \E V \bigl(\widehat{m}_j-m_j\bigr)\leq \bigl(\frac{\eps^2}{a}c_3\left(1+L^{1/2}+c_7\right)+\alpha^{1/2}\bigr)\E V\bigl(\widehat{m}_{j-1}-m_{j-1} \bigr)+c_{11}\eps,
        \end{align}
        where
        $c_{11}=c_5+\frac{\eps}{a}c_3\left(c_9+c_{10}a^{1/2}+c_6+c_8a^{1/2}\right)$.
        Taking $a$ to be large enough to guarantee that $\frac{\eps^2}{a}c_3\left(1+L^{1/2}+c_7\right)\leq \frac{1-\alpha^{1/2}}{2}$, we obtain 
        \begin{align*}
            \E V\bigl(\widehat{m}_j-m_j \bigr)\leq \alpha_* \E V\bigl(\widehat{m}_{j-1}-m_{j-1} \bigr) +c_{11}\eps.
        \end{align*}
        Applying the discrete Gronwall inequality and using the equivalence of $V(\cdot)$ and $\|\cdot\|$ yields the desired result.
\end{proof}
\end{theorem}

The following lemma was used in the proof of Theorem \ref{theorem:general mean field filter accuracy}. This lemma is a generalization of Corollary 4 in \cite{mourtada2022exact} to the setting where the samples are independent, but have different means and thus are not identically distributed.
\begin{lemma}\label{lemma:expected inverse eigenvalue}
    Let $\{u^{(n)}\}_{n=1}^N\in \H$ with empirical covariance satisfying $\Tr(\widehat{C})\leq c\eps^2$. Consider $v^{(n)}=\Psi(P_{\B_V}u^{(n)})-\frac{1}{N}\sum_{n=1}^N\Psi(P_{\B_V}u^{(n)})+\xi^{(n)}$ with $\xi^{(n)}\sim \Nc(0,aP).$ Then, if $N\geq \min\{6k,12\}$ and $a\geq \frac{10NLc\eps^2}{k}$, the sample covariance $\widehat{\Sigma}=\frac{1}{N}\sum_{n=1}^Nv^{(n)}\otimes v^{(n)}$ satisfies
    \begin{align*}
        \E \left[\max\{1,\lambda_{\min}(H\widehat{\Sigma}H^*)^{-q}\}\right]^{1/q}\leq \frac{2C'}{a},
    \end{align*}
    for any $1\leq q \leq N/12$, where $C'$ is a universal constant.
    \begin{proof}
 If we can show that for $N\geq 6k$, for every $t\in (0,1)$,
        \begin{align}\label{eq:lower tail control}
        \mathbb{P}\left(\lambda_{\min} \Bigl(\frac{1}{a}H\widehat{\Sigma}H^* \Bigr) \leq t\right)\leq (C't)^{N/6},
        \end{align}
        for a universal constant $C'$, then the desired result follows exactly as in \cite{mourtada2022exact}. Showing that \eqref{eq:lower tail control} holds is tantamount to verifying that the conclusion of Theorem 4 in \cite{mourtada2022exact} holds in our modified setting. We begin by showing that each $Hv^{(n)}$ satisfies a small ball condition. We have that for every $\theta\in \R^k$ with $\|\theta\|^2=1$
        \begin{align}
        \begin{split}
            \mathbb{P}\left(|\langle \theta,Hv^{(n)}\rangle| \leq t \right)&\leq \sup_{y\in \R}\mathbb{P}\left(|\langle \theta',a^{-1/2}H\xi^{(n)}\rangle -y|\leq t \right)\\
            &\overset{\text{(i)}}{ \leq}\sup_{y\in \R}\int \mathbbm{1}\left\{|z-y|\leq t\right\}\frac{1}{\sqrt{2\pi}}\exp\left(-\frac{1}{2}z^2\right)dz\\
            & \leq\sup_{y\in \R}\int \mathbbm{1}\{|z-a|\leq t\}\frac{1}{\sqrt{2\pi}}dz
            \leq \frac{\sqrt{2}}{\pi}t,
            \end{split}
        \end{align}
        where (i) uses the fact that  $z=\langle\theta,a^{-1/2}H\xi^{(n)}\rangle \sim \Nc(0,1)$ for $\|\theta\|^2=1.$ Given this, we now argue that the proof of Theorem 4 in \cite{mourtada2022exact} goes through essentially as written. As in \cite{mourtada2022exact}, we replace $a^{-1/2}Hv^{(n)}$ by the truncated vectors $\omega^{(n)}:=\max\{1,\frac{\sqrt{k}}{\|a^{-1/2}Hv^{(n)}\|}\}a^{-1/2}Hv^{(n)}$, and define $\widehat{\Sigma}_\omega:=\frac{1}{N}\sum_{n=1}^N\omega^{(n)}(\omega^{(n)})^T$. Since $\widehat{\Sigma}_{\omega}\preceq \widehat{\Sigma},$ we can proceed by establishing a lower bound for $\lambda_{\min }(\widehat{\Sigma}_{\omega})$. For every $\|\theta\|=1$,  $t\in (0,\frac{\pi}{\sqrt{2}})$, and $b\geq 1$,
        \begin{align*}
            \mathbb{P}\left(|\langle \omega^{(n)},\theta\rangle| \leq t \right) & \leq \mathbb{P}\left( |\langle a^{-1/2}Hv^{(n)},\theta\rangle|\leq bt \right)+\mathbb{P}\left(\frac{\sqrt{k}}{\|a^{-1/2}Hv^{(n)}\|}\leq \frac{1}{b} \right) \\
            & \leq \frac{\sqrt{2}}{\pi}bt + \mathbb{P}\left(\|a^{-1/2}Hv^{(n)}\|\geq b \sqrt{k} \right)\\
            & \overset{\text{(i)}}{\leq} \frac{\sqrt{2}}{\pi}bt+\frac{\E \|a^{-1/2}Hv^{(n)}\|^2}{b^2k}\\
            & \leq 
            \frac{\sqrt{2}}{\pi}bt+\frac{k+a^{-1}\|\Psi(P_{\B_V}u^{(n)})-\frac{1}{N}\sum_{n=1}^N\Psi(P_{\B_V}u^{(n)})\|^2}{b^2k}\\
            & \leq  \frac{\sqrt{2}}{\pi}bt+\frac{k+a^{-1}NL\Tr(\widehat{C})}{b^2k}\\
            & \leq  \frac{\sqrt{2}}{\pi}bt+\frac{k+a^{-1}NLc\eps^2}{b^2k}
             \overset{\text{(ii)}}{\leq}\frac{\sqrt{2}}{\pi}bt+\frac{1+\frac{1}{10}}{b^2},
        \end{align*}
        where (i) uses Markov's inequality and (ii) uses our assumption on $a$. Taking $b= \Bigl(\frac{\sqrt{2}}{\pi}t \Bigr)^{-1/3},$ we have
        \begin{align}
            \mathbb{P}\left(|\langle \omega^{(n)},\theta\rangle| \leq t \right)\leq \left(2+\frac{1}{10} \right) \left(\frac{\sqrt{2}}{\pi} t\right)^{2/3},
        \end{align}
        for each $1\leq n \leq N$. Note that this agrees with Equation 69 in \cite{mourtada2022exact} with a constant factor $2+\frac{1}{10}$ in place of the constant factor $2$. Lemma 7 in \cite{mourtada2022exact} then implies that for every $\lambda>0$,
        \begin{align*}
            \E \left[\exp \left(-\lambda \langle \omega^{(n)},\theta\rangle^2 \right)\right]\leq 3 \Bigl(\frac{2}{\pi^2}/\lambda \Bigr)^{1/3},
        \end{align*}
        for each $1\leq n \leq N$. Defining 
        \begin{align*}
            Z_n(\theta)=-\lambda \langle \omega^{(n)},\theta\rangle^2+\frac{1}{3}\log \left(\frac{\lambda}{C^2} \right)-\log\left(2+\frac{1}{10} \right),
        \end{align*}
        for some fixed $\lambda>0$, we necessarily have that $\E \bigl[\exp\bigl(Z_n(\theta)\bigr)\bigr]\leq 1.$ Consequently, letting  
        \begin{align*}
            Z(\theta):=\sum_{n=1}^N Z_n(\theta)=N\left[-\lambda \langle \widehat{\Sigma}_{\omega}\theta,\theta \rangle +\frac{1}{3}\log \left(\frac{\lambda}{C^2
            } \right)-\log \left(2+\frac{1}{10} \right) \right],
        \end{align*}
        the fact that $\omega^{(n)}$ are independent given $\{u^{(n)}\}_{n=1}^N$ implies that for every $\|\theta\|=1$,
        \begin{align}
            \E[Z(\theta)]\leq 1.
        \end{align}
        The PAC-Bayesian argument only considers $Z(\theta)$, and goes through exactly as in \cite{mourtada2022exact}. In the conclusion of the proof, \cite{mourtada2022exact} uses the assumption that $\frac{k}{N}\leq \frac{1}{6}$ to argue that $\log 2 +\frac{k}{N}\log \left( \frac{5N}{2k}\right)+\frac{k}{2N}\leq c_0=1.3.$ Since it also holds that $\log \left(2+\frac{1}{10} \right) +\frac{k}{N}\log \left( \frac{5N}{2k}\right)+\frac{k}{2N}\leq c_0=1.3$, the proof concludes as the proof of Theorem 4 in \cite{mourtada2022exact} to show that \eqref{eq:lower tail control} holds, and our desired result follows.
    \end{proof}
\end{lemma}

\section{Filter accuracy with surrogate models via mean-field theory}\label{sec:proofmain2}

In this section, we prove the long-time accuracy of an ensemble Kalman filters which employs a surrogate dynamics model satisfying Assumption \ref{assumption:surrogate model assumption}. We begin by proving an analogue of Lemma \ref{lemma:covariance trace bound} for the covariance operators generated by Algorithms \ref{algorithm:mean field} and  \ref{algorithm:mean field surrogate}.

\begin{algorithm}[H]
\caption{Gaussian Projected Filter with a Surrogate Model\label{algorithm:mean field surrogate}}
\begin{algorithmic}[1]
\STATE {\bf Initialization}:  $m_0, C_0.$ 
\STATE For $j = 1, 2, \ldots$ do the following:
\STATE {\bf \index{prediction}Prediction}: 
\vspace{-0.3cm}
\begin{align}
\begin{split}
\mu_{j}^s &= \mathbb{E}_{u\sim \Nc(m_{j-1},C_{j-1})}\bigl[\Psi^s(P_{\B_V} u)\bigr],  \\
\Sigma_{j}^s &= \mathbb{E}_{u\sim \Nc(m_{j-1},C_{j-1})}\bigl[\left(\Psi^s(P_{\B_V}u)-\mu_{j} \right)\otimes\left(\Psi^s(P_{\B_V}u)-\mu_{j} \right)\bigr].
\end{split}
\end{align}
\vspace{-0.5cm}
\STATE{{{\bf \index{analysis}Analysis}}}: 
\vspace{-0.3cm}
	\begin{align}
	\begin{split}
        m_{j}^s&=\mu_{j}^s+\msK(\Sigma_{j}^s + Q)\left(y_{j} -H\mu_{j}^s\right),\\
        C_{j}^s & =  \bigl( I - \msK(\Sigma_{j}^s + Q) H \bigr) (\Sigma_{j}^s + Q).
\end{split}
\end{align}
\vspace{-0.4cm}
\STATE {\bf Output}: Means and covariances  $\{m_j^s, C_j^s\}_{j=1}^\infty.$
\end{algorithmic}
\end{algorithm}

The next result, analogous to Lemma \ref{lemma:covariance trace bound}, shows that the trace of the covariances generated by Algorithms \ref{algorithm:finite ensemble surrogate} and \ref{algorithm:mean field surrogate} are controlled at the noise level $\eps$ and the surrogate error $\delta$ in the long-time asymptotic. 

\begin{lemma}\label{lemma:surrogate covariance trace bound}
    Under Assumptions \ref{assumption:ball and squeezing} and \ref{assumption:surrogate model assumption}, the covariance operators output by Algorithm \ref{algorithm:mean field surrogate} satisfy
    \begin{align}\label{eq:surrogate mean field trace bound}
        \Tr\left(C^s_j \right)\leq \alpha_*^j\Tr \left(C_0 \right)+\left(1+\alpha_*\beta \right)\frac{1-\alpha_*^j}{1-\alpha_*}\eps^2\Tr\left(R \right)+ 4\frac{1+\alpha}{1-\alpha}\frac{1-\alpha_*^j}{1-\alpha_*}\delta^2,
    \end{align}
    and those output by Algorithm \ref{algorithm:finite ensemble surrogate} satisfy
    \begin{align}\label{eq:surrogate finite ensemble trace bound}
        \Tr \bigl(\widehat{C}^s_j \bigr)\leq \alpha_*^j\Tr \bigl(\widehat{C}_0 \bigr)+\left(1+\alpha_*\beta \right)\frac{1-\alpha_*^j}{1-\alpha_*}\eps^2\Tr\left(R \right)+4\frac{1+\alpha}{1-\alpha}\frac{1-\alpha_*^j}{1-\alpha_*}\delta^2,
    \end{align}
    where $\alpha_*:=\frac{1+\alpha}{2}.$
    \begin{proof}
        Since $P$ is an orthogonal projection and $Q = aP,$ we have
            \begin{align}\label{eq:auxbound surrogate}
            \begin{split}
                    \Tr(C_{j}^s)&=\Tr(P C_{j}^sP^*)+\Tr \bigl((I-P)C_{j}^s(I-P)^* \bigr)\\
            & \le \Tr(P C_{j}^sP^*)+\Tr \bigl((I-P)\Sigma_j^s (I-P)^* \bigr).
            \end{split}
            \end{align}
We next bound in turn the two terms in the right-hand side. The first term proceeds exactly as in Lemma \ref{lemma:covariance trace bound}. Let $\{\phi_i\}_{i=1}^k$ be an orthonormal basis for the range of $P$, which we can extend to a basis for $\H$ via Gram-Schmidt by the separability of $\H$. We then compute that
            \begin{align}
            \begin{split}\label{eq:firstterm surrogate}
\Tr\left(PC_{j}^sP^*\right)&=\sum_{i=1}^{\infty}\langle \phi_i, PC_{j}^sP^* \phi_i \rangle
                =\sum_{i=1}^{k}\langle \phi_i, PC_{j}^sP^* \phi_i \rangle\\
                &=\sum_{i=1}^k\langle H \phi_i, HC_j^sH^* H\phi_i \rangle_{\R^k}= \Tr\left(HC_j^sH^* \right) \le \eps^2\Tr\left(R\right),
            \end{split}
            \end{align}
           where we have used that $\{H\phi_i\}_{i=1
            }^k$ is an orthonormal basis for $\R^k$ and that $  \eps^2 R - HC_{j}^sH^*$ is positive semi-definite in the last inequality.

             For the second term in the right-hand side of \eqref{eq:auxbound surrogate},   note that

        \begin{align}\label{eq:secondterm surrogate}
            \begin{split}
                \Tr \Bigl( (I-P)\Sigma_j^s(I-P)^* \Bigr)&=\mathbb{E}_{\Nc(m_{j-1}^s,C_{j-1}^s)} \Bigl[ \bigl\| (I-P) \bigl(\Psi^s(P_{\B_V} u)-\mu_j \bigr) \bigr\|^2 \Bigr]\\
                & \hspace{-2cm} \leq \mathbb{E}_{\Nc(m_{j-1}^s,C_{j-1}^s)}\Bigl[ \bigl\|(I-P)\bigl(\Psi^s(P_{\B_V} u)-\Psi(P_{\B_V}m^s_{j-1}) \bigr) \bigr\|^2\Bigr]\\
                & \hspace{-2cm} \leq \mathbb{E}_{\Nc(m_{j-1}^s,C_{j-1}^s)}\Bigl[ \bigl\|(I-P)\bigl(\Psi^s(P_{\B_V}u)-\Psi(P_{\B_V}u) 
 \bigr) \bigr\|^2\\
                &\hspace{-2cm}  +2\bigl\|(I-P)\bigl(\Psi^s(P_{\B_V}u)-\Psi(P_{\B_V}u) \bigr)\bigr\|\bigl\|(I-P)(\Psi(P_{\B_V}u)-\Psi(P_{\B_V}m^s_{j-1}) \bigr\|\\
                &\hspace{-2cm}  + \bigl\|(I-P)\bigl(\Psi(P_{\B_V}u)-\Psi(P_{\B_V}m^s_{j-1} \bigr) \bigr\|^2 \Bigr]\\
                &\hspace{-2cm}  \overset{\text{(i)}}{\leq}  \mathbb{E}_{\Nc(m_{j-1}^s,C_{j-1}^s)}\biggl[4\frac{1+\alpha}{1-\alpha} \|(I-P) \bigl(\Psi^s(P_{\B_v}u)-\Psi(P_{\B_V}u) \bigr) \bigr\|^2\\
                &\hspace{-2cm} +\frac{1+\alpha}{2\alpha}\bigl\|(I-P) \bigl(\Psi(P_{\B_V}u)-\Psi(P_{\B_V}m^s_{j-1}) \bigr) \bigr\|^2 \biggr]\\
                &\hspace{-2cm}  \overset{\text{(ii)}}{\leq} 4\frac{1+\alpha}{1-\alpha}\delta^2+\mathbb{E}_{\Nc(m_{j-1}^s,C_{j-1}^s)}\biggl[\frac{1+\alpha}{2\alpha}\bigl\|(I-P) \bigl(\Psi(P_{\B_V}u)-\Psi(P_{\B_V}m^s_{j-1}) \bigr) \bigr\|^2 \biggr]\\
                &\hspace{-2cm}  \overset{\text{(iii)}}{\leq }4\frac{1+\alpha}{1-\alpha}\delta^2+\mathbb{E}_{\Nc(m_{j-1}^s,C_{j-1}^s)}\biggl[\frac{1+\alpha}{2}\Bigl(\bigl\|u-m_{j-1}^s\|^2+\beta\bigl\|(I-P) (u-m_{j-1}^s)\bigr\|^2 \Bigr) \biggr]\\
                & \hspace{-2cm} = 4\frac{1+\alpha}{1-\alpha}\delta^2+\alpha_* \Tr\left(C_{j-1}^s \right)+\alpha_*\beta \Tr\left(PC_{j-1}^sP^* \right),
                \end{split}
            \end{align}
            where (i) uses Young's inequality, (ii) uses Assumption \ref{assumption:surrogate model assumption}, and (iii) uses Assumption \ref{assumption:ball and squeezing}. Combining \eqref{eq:auxbound surrogate}, \eqref{eq:firstterm surrogate}, and \eqref{eq:secondterm surrogate}, we deduce that 
         \begin{align*}
                \Tr(C_{j}^s)& \leq \eps^2\Tr(R)+4\frac{1+\alpha}{1-\alpha}\delta^2+\alpha_* \Tr\left(C_{j-1}^s \right)+\alpha_*\beta \Tr\left(PC_{j-1}^sP^* \right)\\
                & \leq 4\frac{1+\alpha}{1-\alpha}\delta^2+(1+\alpha_* \beta)\eps^2\Tr (R)+\alpha_* \Tr (C^s_{j-1}).
            \end{align*}
            Applying the discrete Gronwall lemma yields \eqref{eq:surrogate mean field trace bound}. The proof of \eqref{eq:surrogate finite ensemble trace bound} is similar.
    \end{proof}
\end{lemma}

The next result establishes long-time accuracy for the idealized mean-field Algorithm \ref{algorithm:mean field surrogate}.
\begin{theorem}\label{thm:surrogate1}
    Suppose that Assumptions \ref{assumption:ball and squeezing} and \ref{assumption:surrogate model assumption} hold. Then, if the inflation parameter $a>0$ is sufficiently large, there is a constant $\mathsf{C}_3$ independent of $\eps$ and $\delta$ such that
    \begin{align}
        \limsup_{j\to \infty}\E \|m_j-m_j^s\|\leq \mathsf{C}_3(\eps+\delta).
    \end{align}
    \begin{proof}
        Assume $a$ is taken to be at least as large as in Theorem \ref{theorem:general mean field filter accuracy}. From Lemma \ref{lemma:covariance trace bound}, Lemma \ref{lemma:surrogate covariance trace bound} and Theorem \ref{theorem:general mean field filter accuracy} there exists a time $j_*$ such that for all $j\geq j_*$, $\Tr(C_j)\leq c_1\eps^2$ with $c_1:=1+\frac{1+\alpha\beta}{1-\alpha}\Tr(R),$ $\Tr(C_j^s)\leq c_2(\eps^2+\delta^2)$ with $c_2=2\max\{(1+\alpha_*\beta)\frac{1}{1-\alpha_*}\Tr(R),4\frac{1+\alpha}{1-\alpha}\frac{1}{1-\alpha_*}\},$ and that $\E V(m_j-u_j)\leq 2\mathsf{C}_1\eps.$ For $j>j^*,$ we consider the difference
        \begin{align*}
            m_j^s-m_j&=\left(I-\mathscr{K}(\Sigma_j^s+Q)H \right)\mu_j^s+\mathscr{K}(\Sigma_j^s+Q)y_j-\left(I-\mathscr{K}(\Sigma_j+Q) \right)\mu_j-\mathscr{K}(\Sigma_j+Q)y_j\\
            & \hspace{-1cm} = \Bigl(I -\mathscr{K}(\Sigma_j^s+Q)H \Bigr)(\mu_j^s-\mu_j)+ \Bigl(I -\mathscr{K}(\Sigma_j^s+Q)H-I +\mathscr{K}(\Sigma_j+Q)H \Bigr)(\mu_j-H^*y_j ).
        \end{align*}
        Taking the $V$ norm of both sides yields
        \begin{align}\label{eq:mf surrogate bound 1}
        \begin{split}
            V\bigl(m_j^s-m_j\bigr)&\leq V\Bigl(\bigl(I -\mathscr{K}(\Sigma_j^s+Q)H \bigr)\bigl(\mu_j^s-\mu_j \bigr) \Bigr)\\
            & +V \Bigl( \bigl(I -\mathscr{K}(\Sigma_j^s+Q)H-I +\mathscr{K}(\Sigma_j+Q)H \bigr) \bigl(\mu_j-H^*y_j \bigr)\Bigr).
            \end{split}
        \end{align}
        The first of the terms on the right-hand side of \eqref{eq:mf surrogate bound 1} can be bounded as
        \begin{align*}
            V\Bigl(\bigl(I -\mathscr{K}(\Sigma_j^s+Q)H \bigr)\bigl(\mu_j^s-\mu_j\bigr) \Bigr) &\leq V \Bigl(\bigl(I -\mathscr{K}(\Sigma_j^s+Q)H \bigr) \bigl(\mu_j^s-\Psi^s(P_{\B_V}m_{j-1}^s) \bigr) \Bigr)\\
            &\hspace{-1cm}  +V \Bigl(\bigl(I -\mathscr{K}(\Sigma_j^s+Q)H \bigr)\bigl(\Psi^s(P_{\B_V}m_{j-1}^s)-\Psi(P_{\B_V}m_{j-1}^s) \bigr) \Bigr)\\
            &\hspace{-1cm}  + V \Bigl(\bigl(I -\mathscr{K}(\Sigma_j^s+Q)H \bigr)\bigl(\Psi(P_{\B_V}m_{j-1}^s)-\Psi(P_{\B_V}m_{j-1}) \bigr) \Bigr)\\
            &\hspace{-1cm}  + V \Bigl(\bigl(I -\mathscr{K}(\Sigma_j^s+Q)H \bigr)\bigl(\Psi(P_{\B_V}m_{j-1})-\mu_j \bigr) \Bigr).
        \end{align*}
        We bound each of these terms in turn. By an analogous argument to that in Theorem \ref{theorem:general mean field filter accuracy},  
        \begin{align}\label{eq:mf P minus Kalman gain surrogate}
            \|P -\mathscr{K}(\Sigma_j^s+Q)H\|_{V,V}\leq \frac{1}{a}\Bigl(\eps^2\|R\|_{op}+L_sc_1(\eps^2+\delta^2) \Bigr)
        \end{align}
        and 
        \begin{align}\label{eq:mf I minus Kalman gain surrogate}
            \|I -\mathscr{K}(\Sigma_j^s+Q)H\|_{V,V}\leq 1+\frac{1}{a}\Bigl(\eps^2\|R\|_{op}+L_sc_1(\eps^2+\delta^2) \Bigr).
        \end{align}
       Defining $c_3:=\Bigl(\eps^2\|R\|_{op}+L_sc_2(\eps^2+\delta^2) \Bigr)$, we have
        \begin{align}\label{eq:mf surrogate bound 2}
        \begin{split}
            V \Bigl(\bigl(I -\mathscr{K}(\Sigma_j^s+Q)H \bigr)\bigl(\mu_j^s-\Psi^s(m_{j-1}^s)\bigr) \Bigr)&\leq \|I -\mathscr{K}(\Sigma_j^s+Q)H\|_{V,V} V\Bigl(\mu_j^s-\Psi^s(m_{j-1}^s) \Bigr)\\ 
            & \leq \Bigl(1+\frac{c_3}{a}\Bigr)V\Bigl(\mu_j^s-\Psi^s(m_{j-1}^s) \Bigr) \\
            & \leq \Bigl(1+\frac{c_3}{a} \Bigr)L_s^{1/2}\Bigl(\Tr(C_{j-1}^s)+\Tr(HC_{j-1}H^*)\Bigr)^{1/2}\\
            & \leq \Bigl(1+\frac{c_3}{a} \Bigr)L_s^{1/2}\Bigl(c_2(\eps^2+\delta^2)+\eps^2\Tr(R) \Bigr)^{1/2}\\
            &\leq c_4(\eps+\delta),
            \end{split}
        \end{align}
        with $c_4:=(1+\frac{c_3}{a})L_s^{1/2}(c_2+\Tr(R))^{1/2}.$ Next, we have that
        \begin{align}\label{eq:mf surrogate bound 3}
            \begin{split}
                &V \Bigl(\bigl(I -\mathscr{K}(\Sigma_j^s+Q)H \bigr)\bigl(\Psi^s(P_{\B_V}m_{j-1}^s)-\Psi(P_{\B_V}m_{j-1}^s) \bigr) \Bigr) \\
                &\leq V \Bigl(\bigl(I -\mathscr{K}(\Sigma_j^s+Q)H \bigr)P \bigl(\Psi^s(P_{\B_V}m_{j-1}^s)-\Psi(P_{\B_V}m_{j-1}^s) \bigr) \Bigr)\\
                & + V \Bigl(\bigl(I -\mathscr{K}(\Sigma_j^s+Q)H \bigr)(I-P) \bigl(\Psi^s(P_{\B_V}m_{j-1}^s)-\Psi(P_{\B_V}m_{j-1}^s) \bigr) \Bigr)\\
                & \leq \|P-\mathscr{K}(\Sigma_j^s+Q)H\|_{V,V}V\Bigl(P \bigl(\Psi^s(P_{\B_V}m_{j-1}^s)-\Psi(P_{\B_V}m_{j-1}^s) \bigr)\Bigr)\\
                & + V\Bigl((I-P) \bigl(\Psi^s(P_{\B_V}m_{j-1}^s)-\Psi(P_{\B_V}m_{j-1}^s) \bigr) \Bigr)\\
                & \overset{\text{(i)}}{\leq} \|P-\mathscr{K}(\Sigma_j^s+Q)H\|_{V,V}\sqrt{2}\kappa +\delta\\
                & \leq \frac{\eps^2\|R\|_{op}+L_sc_1(\eps^2+\delta^2)}{a}\sqrt{2}\kappa+\delta
                 \leq c_5(\eps+\delta),
            \end{split}
        \end{align}
        where (i) uses Assumption \ref{assumption:surrogate model assumption} and $c_5:=2\sqrt{2}\kappa\frac{\|R\|_{op}+L_sc_1}{a}(\eps+\delta)+1.$ Similarly, we compute 
        \begin{align}\label{eq:mf surrogate bound 4}
            \begin{split}
                &V \Bigl(\bigl(I -\mathscr{K}(\Sigma_j^s+Q)H \bigr) \bigl(\Psi(P_{\B_V}m_{j-1}^s)-\Psi(P_{\B_V}m_{j-1}) \bigr) \Bigr) \\
                & \leq V \Bigl(\bigl(I -\mathscr{K}(\Sigma_j^s+Q)H \bigr)P \bigl(\Psi(P_{\B_V}m_{j-1}^s)-\Psi(P_{\B_V}m_{j-1}) \bigr)\Bigr)\\
                &+V \Bigl(\bigl(I -\mathscr{K}(\Sigma_j^s+Q)H \bigr)(I-P) \bigl(\Psi(P_{\B_V}m_{j-1}^s)-\Psi(P_{\B_V}m_{j-1}) \bigr)\Bigr)\\
                & \leq \|P-\mathscr{K}(\Sigma_j^s+Q)H\|_{V,V}V\Bigl(\Psi(P_{\B_V}m_{j-1}^s)-\Psi(P_{\B_V}m_{j-1})\Bigr)\\
                & + V\Bigl((I-P) \bigl(\Psi(P_{\B_V}m_{j-1}^s)-\Psi(P_{\B_V}m_{j-1}) \bigr)\Bigr)\\
                & \overset{\text{(i)}}{\leq} \|P-\mathscr{K}(\Sigma_j^s+Q)H\|_{V,V}L^{1/2}V\Bigl(m_{j-1}^s-m_{j-1}\Bigr)+\alpha^{1/2}V\bigl(m_{j-1}^s-m_{j-1}\bigr)\\
                & \leq \Bigl(\frac{\eps^2\|R\|_{op}+L_sc_2(\eps^2+\delta^2)}{a}+\alpha^{1/2}\Bigr)V\bigl(m_{j-1}^s-m_{j-1}\bigr)\\
                & =\Bigl(\frac{c_3}{a}+\alpha^{1/2}\Bigr)V\bigl(m_{j-1}^s-m_{j-1}\bigr),
            \end{split}
        \end{align}
where (i) uses Assumption \ref{assumption:ball and squeezing}. Finally, we have that 
\begin{align}\label{eq:mf surrogate bound 5}
    \begin{split}
        V \Bigl(\bigl(I -\mathscr{K}(\Sigma_j^s+Q)H \bigr)(\Psi(P_{\B_V}m_{j-1})-\mu_j) \Bigr)& \leq \Bigl(1+\frac{c_2}{a}\Bigr)V\Bigl(\Psi(P_{\B_V}m_{j-1})-\mu_j\Bigr)\\
        &  \hspace{-1cm} \leq \Bigl(1+\frac{c_2}{a}\Bigr)L^{1/2}\Bigl(\Tr(C_{j-1})+\Tr(HC_{j-1}H^*)\Bigr)^{1/2}\\
        & \hspace{-1cm} \leq \Bigl(1+\frac{c_2}{a}\Bigr)L^{1/2}\Bigl(c_1\eps^2+\Tr(R)\eps^2\Bigr)^{1/2}  = c_6\eps,
    \end{split}
\end{align}
        with $c_6:=(1+\frac{c_2}{a})L^{1/2}(c_1+\Tr(R))^{1/2}.$ Combining \eqref{eq:mf surrogate bound 2}, \eqref{eq:mf surrogate bound 3}, \eqref{eq:mf surrogate bound 4}, and \eqref{eq:mf surrogate bound 5}, we have
        \begin{align}\label{eq:mf surrogate bound 6}
             V\Bigl(\bigl(I -\mathscr{K}(\Sigma_j^s+Q)H \bigr)(\mu_j^s-\mu_j) \Bigr)\leq \Bigl(\frac{c_3}{a}+\alpha^{1/2} \Bigr)V\bigl(m^s_{j-1}-m_{j-1} \bigr)+c_7(\eps+\delta),
        \end{align}
        where $c_7:=c_4+c_5+c_6.$ We now bound the second term in \eqref{eq:mf surrogate bound 1}.
        \begin{align*}
            \begin{split}
                &V \Bigl( \bigl(I -\mathscr{K}(\Sigma_j^s+Q)H-I +\mathscr{K}(\Sigma_j+Q)H \bigr) \bigl(\mu_j-H^*y_j \bigr)\Bigr)\\
                & =  V \Bigl( \bigl(P -\mathscr{K}(\Sigma_j^s+Q)H-P +\mathscr{K}(\Sigma_j^s+Q)H \bigr)\bigl(P\mu_j-H^*y_j \bigr)\Bigr)\\
                & \leq \Bigl(\|P -\mathscr{K}(\Sigma_j^s+Q)H\|_{V,V}+\|P -\mathscr{K}(\Sigma_j+Q)H\|_{V,V}\Bigr) V\bigl(P\mu_j-H^*y_j\bigr)\\
                & \leq \biggl( \frac{\eps^2\|R\|_{op}+L_sc_1(\eps^2+\delta^2)}{a}+\frac{\eps^2(\|R\|_{op}+L(1+\frac{1+\alpha \beta}{1-\alpha}\Tr(R)))}{a}\biggr)V\bigl(P\mu_j-H^*y_j\bigr)\\
                & \overset{\text{(i)}}{\leq}  c_8(\eps+\delta)V\bigl(P\mu_j-H^*y_j\bigr)\\
                & \leq c_8(\eps+\delta)\Bigl[V\Bigl(\mu_j-\Psi(P_{\B_V}m_{j-1})\Bigr)+V\Bigl(\Psi(P_{B_V}m_{j-1})-\Psi(u_{j-1})\Bigr)+ V(H^*\eta) \Bigr]\\
                &\overset{\text{(ii)}}{\leq} c_8(\eps+\delta)\Bigl[L^{1/2}\Bigl(\Tr(c_1\eps^2+\Tr(R)\eps^2\Bigr)^{1/2}+V\Bigl(\Psi(P_{B_V}m_{j-1})-\Psi(u_{j-1})\Bigr)+ V(H^*\eta)\Bigr],            \end{split}
        \end{align*}
        where in (i) we define $c_8:=\frac{2\eps\|R\|_{op}+2L_sc_1(\eps+\delta)+L(1+\frac{1+\alpha \beta}{1-\alpha}\Tr(R))}{a}$, and (ii) uses the same argument as in \eqref{eq:mf surrogate bound 5}. Taking expectations, we have
        \begin{align}\label{eq:mf surrogate bound 7}
            \begin{split}
                &\E V \Bigl( \bigl(I -\mathscr{K}(\Sigma_j^s+Q)H-I +\mathscr{K}(\Sigma_j+Q)H \bigr) \bigl(\mu_j-H^*y_j \bigr)\Bigr) \\
                & \leq c_8(\eps+\delta)\Bigl(L^{1/2}\Bigl(\Tr(c_1\eps^2+\Tr(R)\eps^2\Bigr)^{1/2}+L^{1/2}2\mathsf{C}_1\eps +\eps \sqrt{2\Tr(R)}\Bigr)\\
                & = c_9(\eps+\delta)
            \end{split}
        \end{align}
        with $c_9:=c_8\Bigl(L^{1/2}\Bigl(\Tr(c_1\eps^2+\Tr(R)\eps^2\Bigr)^{1/2}+L^{1/2}2\mathsf{C}_1\eps +\eps \sqrt{2\Tr(R)}\Bigr).$ Combining \eqref{eq:mf surrogate bound 6} and \eqref{eq:mf surrogate bound 7},
        \begin{align*}
            \E V\bigl(m_j^s-m_j\bigr)\leq \Bigl(\frac{c_3}{a}+\alpha^{1/2} \Bigr)V\bigl(m^s_{j-1}-m_{j-1} \bigr)+(c_7+c_9)(\eps+\delta),
        \end{align*}
        and we see that taking $a\geq \frac{2c_3}{1-\alpha^{1/2}},$
        \begin{align*}
            \E V\bigl(m_j^s-m_j\bigr)\leq \alpha_*V\bigl(m^s_{j-1}-m_{j-1} \bigr)+(c_7+c_9)(\eps+\delta)
        \end{align*}
        with $\alpha_*=\frac{1+\alpha^{1/2}}{2}.$ The discrete-time Gronwall lemma and equivalence of $V(\cdot)$ and $\|\cdot\|$ completes the proof.
    \end{proof}
\end{theorem}

We now show that the means output by Algorithm \ref{algorithm:finite ensemble surrogate} and those output by the idealized mean-field Algorithm \ref{algorithm:mean field surrogate} remain close in the long-time asymptotic. The proof of our second main result, Theorem \ref{th:main2}, follows directly by combining Theorem \ref{thm:surrogate1} with Theorem \ref{theorem:general mean field filter accuracy}
and the next result: 

\begin{theorem}
    Suppose that Assumptions \ref{assumption:ball and squeezing} and \ref{assumption:surrogate model assumption} hold and that $N\geq 6k$. Then, if the inflation parameter $a>0$ is sufficiently large, there is a constant $\mathsf{C}_4$ independent of $\eps$ and $\delta$ such that
    \begin{align}
        \limsup_{j\to \infty}\E \|\widehat{m}^s_j-m^s_j\|\leq \mathsf{C}_4(\eps+\delta).
    \end{align}
    \begin{proof}
        Assume that $a$ is taken to be at least as large as in Theorem \ref{thm:surrogate1}. From Lemma \ref{lemma:surrogate covariance trace bound} and Theorem \ref{thm:surrogate1} there exists a time $j_*$ such that for $j\geq j_*$, $\Tr(C^s_j)\leq c_1(\eps^2+\delta^2)$ and $\Tr(\widehat{C}_j^s)\leq c_1(\eps^2+\delta^2)$ with $c_1:=2\max\{\frac{1+\alpha_*\beta}{1-\alpha_*}\Tr(R),4\frac{1+\alpha}{1-\alpha}\frac{1}{1-\alpha^*}\},$ and that $\E V (m^s_j-u_j)\leq 2(\mathsf{C}_1+\mathsf{C}_3)(\eps+\delta)$. For $j>j^*,$ we consider the difference
        \begin{align*}
            \widehat{m}_j^s-m^s_j&=\bigl(I-\mathscr{K}(\widehat{\Sigma}_j^s)H \bigr)\widehat{\mu}_j^s+\mathscr{K}(\widehat{\Sigma}_j^s)y_j-\bigl(I-\mathscr{K}(\Sigma^s_j) \bigr)\mu^s_j-\mathscr{K}(\Sigma^s_j)y_j\\
            & = \bigl(I -\mathscr{K}(\widehat{\Sigma}_j^s)H \bigr)(\widehat{\mu}_j^s-\mu^s_j)+ \bigl(I -\mathscr{K}(\widehat{\Sigma}_j^s)H-I +\mathscr{K}(\Sigma^s_j)H \bigr)(\mu^s_j-H^*y_j ).
        \end{align*}
        Taking the $V$ norm and expectation of both sides gives
        \begin{align}\label{eq:fp surrogate bound 1}
            \begin{split}
                \E V\bigl(\widehat{m}_j^s-m^s_j\bigr) & \leq \E V\Bigl( \bigl(I-\mathscr{K}(\widehat{\Sigma}_j^s)H \bigr) \bigl(\widehat{\mu}_j^s-\mu^s_j \bigr)\Bigr) \\
                & + \E V\Bigl(\bigl(I -\mathscr{K}(\widehat{\Sigma}_j^s)H-I +\mathscr{K}(\Sigma^s_j)H \bigr) \bigl(\mu^s_j-H^*y_j \bigr)\Bigr).
            \end{split}
        \end{align}
        The first term on the right-hand side of \eqref{eq:fp surrogate bound 1} can be bounded as
    \begin{align*}
        &\E V\biggl( \bigl(I-\mathscr{K}(\widehat{\Sigma}_j^s)H \bigr)\bigl(\widehat{\mu}_j^s-\mu^s_j \bigr)\biggr)\\ & \leq \E V\biggl( \bigl(I-\mathscr{K}(\widehat{\Sigma}_j^s)H \bigr) \Bigl(\frac{1}{N}\sum_{n=1}^N\Psi^s(P_{\B_V}u^{(n),s}_{j-1})-\Psi^s(P_{\B_V}\widehat{m}_{j-1}^s)\Bigr)\biggr) \\
        & + \E V\biggl( \bigl(I-\mathscr{K}(\widehat{\Sigma}_j^s)H \bigr) \Bigl(\Psi^s(P_{\B_V}\widehat{m}^s_{j-1})-\Psi(P_{\B_V}\widehat{m}^s_{j-1}) \Bigr)\biggr)\\
        & + \E V\biggl( \bigl(I-\mathscr{K}(\widehat{\Sigma}_j^s)H \bigr) \Bigl(\Psi(P_{\B_V}\widehat{m}_{j-1}^s)-\Psi(P_{\B_V}m^s_{j-1})\Bigr)\biggr)\\
        & +\E V\biggl( \bigl(I-\mathscr{K}(\widehat{\Sigma}_j^s)H \bigr) \Bigl(\Psi(P_{\B_V}m^s_{j-1})-\Psi^s(P_{\B_V}m^s_{j-1}) \Bigr)\biggr)\\
        & + \E V\biggl( \bigl(I-\mathscr{K}(\widehat{\Sigma}_j^s)H \bigr) \Bigl(\Psi^s(P_{\B_V}m^s_{j-1})-\mu^s_{j} \Bigr)\biggr)\\
        & +\E V\biggl( \bigl(I-\mathscr{K}(\widehat{\Sigma}_j^s)H \bigr)\frac{1}{N} \sum_{n=1}^N \xi^{(n)}\biggr).\\
    \end{align*}
    We bound each of these terms in turn. By a similar argument to that in the proof of Theorem \ref{theorem:general mean field filter accuracy}, we have
    \begin{align}\label{eq: P minus Kalman gain norm bound}
        \|P-\mathscr{K}(\widehat{\Sigma}_j^s)H\|_{V,V}\leq \frac{\eps^2+\delta^2}{\lambda_{\min}(H\widehat{\Sigma}_j^sH^*)}\Bigl(\|R\|_{op}+L_sc_1\Bigr),
    \end{align}
    and 
    \begin{align}
        \|I-\mathscr{K}(\widehat{\Sigma}_j^s)H\|_{V,V}\leq 1+\frac{\eps^2+\delta^2}{\lambda_{\min}(H\widehat{\Sigma}_j^sH^*)}\Bigl(\|R\|_{op}+L_sc_1\Bigr).
    \end{align}
    Since $\widehat{\Sigma}_j^s$ is conditionally independent of $\frac{1}{N}\sum_{n=1}^N\Psi^s(P_{\B_V}u^{(n),s}_{j-1})-\Psi^s(P_{\B_V}\widehat{m}_{j-1}^s)$ given $\mathsf{U_{j-1}^s}$, 
    \begin{align}
    \begin{split}
        &\E V\biggl( \bigl(I-\mathscr{K}(\widehat{\Sigma}_j^s)H \bigr) \Bigl(\frac{1}{N}\sum_{n=1}^N\Psi^s(P_{\B_V}u^{(n),s}_{j-1})-\Psi^s(P_{\B_V}\widehat{m}_{j-1}^s)\Bigr)\biggr)\\
        & \leq \E \biggl[\mathbb{E}^\mathsf{U_{j-1}^s}\Bigl[1+\frac{\eps^2+\delta^2}{\lambda_{\min}(H\widehat{\Sigma}_j^sH^*)}\Bigl(\|R\|_{op}+L_sc_1\Bigr) \Bigr]\mathbb{E}^\mathsf{U_{j-1}^s} \Bigl[V\Bigl(\frac{1}{N}\sum_{n=1}^N\Psi^s(P_{\B_V}u^{(n),s}_{j-1})-\Psi^s(P_{\B_V}\widehat{m}_{j-1}^s)\Bigr)\Bigr]\biggr]\\
        & \leq \E \biggl[\mathbb{E}^\mathsf{U_{j-1}^s} \Bigl[1+\frac{\eps^2+\delta^2}{\lambda_{\min}(H\widehat{\Sigma}_j^sH^*)}\Bigl(\|R\|_{op}+L_sc_1\Bigr) \Bigr]\mathbb{E}^\mathsf{U_{j-1}^s}\Bigl[L_s^{1/2}\Bigl(\Tr(\widehat{C}_{j-1}^s)+\Tr (H\widehat{C}_{j-1}^sH^*)\Bigr)^{1/2}\Bigr]\biggr]\\
        & \leq \E \biggl[\mathbb{E}^\mathsf{U_{j-1}^s}\Bigl[1+\frac{\eps^2+\delta^2}{\lambda_{\min}(H\widehat{\Sigma}_j^sH^*)}\Bigl(\|R\|_{op}+L_sc_1\Bigr) \Bigr]L_s^{1/2}\sqrt{2c_2}(\eps+\delta)\biggr].
        \end{split}
    \end{align}
    Choosing $a$ to guarantee that $\frac{10 NL_sc_1(\eps^2+\delta^2)}{ak}\leq 1$, Lemma \ref{lemma:expected inverse eigenvalue} gives us that
    \begin{align*}
        \mathbb{E}^\mathsf{U_{j-1}} \Bigl[1+\frac{\eps^2+\delta^2}{\lambda_{\min}(H\widehat{\Sigma}_j^sH^*)}\Bigl(\|R\|_{op}+L_sc_1\Bigr) \Bigr]\leq \biggl(1+\frac{2C'(\eps^2+\delta^2)}{a}\Bigl(\|R\|_{op}+L_sc_1\Bigr) \biggr),
    \end{align*}
    and consequently
    \begin{align}\label{eq:fp surrogate bound 2}
        \begin{split}
            &\E V\biggl( \bigl(I-\mathscr{K}(\widehat{\Sigma}_j^s)H \bigr) \Bigl(\frac{1}{N}\sum_{n=1}^N\Psi^s(P_{\B_V}u^{(n),s}_{j-1})-\Psi^s(P_{\B_V}\widehat{m}_{j-1}^s)\Bigr)\biggr)\\
            & \leq\Bigl(1+\frac{2C'(\eps^2+\delta^2)}{a}\Bigl(\|R\|_{op}+L_sc_2\Bigr) \Bigr)L_s^{1/2}\sqrt{2c_2}(\eps+\delta)
             = c_2(\eps+\delta),
        \end{split}
    \end{align}
    where $c_2:=\biggl(1+\frac{2C'(\eps^2+\delta^2)}{a}\Bigl(\|R\|_{op}+L_sc_1\bigr) \biggr)L_s^{1/2}\sqrt{2c_1}$. Similarly, for the second term we have
    \begin{align}\label{eq:fp surrogate bound 3}
        \begin{split}
            &\E V\biggl( \bigl(I-\mathscr{K}(\widehat{\Sigma}_j^s)H \bigr) \Bigl(\Psi^s(P_{\B_V}\widehat{m}^s_{j-1})-\Psi(P_{\B_V}\widehat{m}^s_{j-1}) \Bigr)\biggr)\\
            &\leq \E V\biggl(\bigl(P-\mathscr{K}(\widehat{\Sigma}_j^s)H \bigr) \Bigl(P\Psi^s(P_{\B_V}\widehat{m}^s_{j-1})-P\Psi(P_{\B_V}\widehat{m}^s_{j-1}) \Bigr) \biggr)\\
            &+\E V\biggl((I-P)\Psi^s(P_{\B_V}\widehat{m}^s_{j-1})-(I-P)\Psi(P_{\B_V}\widehat{m}^s_{j-1})\biggr)\\
            & \leq \E \biggl[\mathbb{E}^\mathsf{U_{j-1}^s} \Bigl[\frac{\eps^2+\delta^2}{\lambda_{\min}(H\widehat{\Sigma}_j^sH^*)}\Bigl(\|R\|_{op}+L_sc_1\Bigr) \Bigr]\mathbb{E}^\mathsf{U_{j-1}^s}\Bigl[\Psi^s(P_{\B_V}\widehat{m}^s_{j-1})-\Psi(P_{\B_V}\widehat{m}^s_{j-1})\Bigr]\biggr]\\
            &+ \E V\Bigl((I-P)\Psi^s(P_{\B_V}\widehat{m}^s_{j-1})-(I-P)\Psi(P_{\B_V}\widehat{m}^s_{j-1})\Bigr)\\
            &\leq \E \biggl[\mathbb{E}^\mathsf{U_{j-1}^s} \Bigl[\frac{\eps^2+\delta^2}{\lambda_{\min}(H\widehat{\Sigma}_j^sH^*)}\Bigl(\|R\|_{op}+L_sc_1\Bigr) \Bigr]\sqrt{2}\kappa \biggr]+\delta\\
            & \leq \frac{2\sqrt{2}\kappa C'(\eps^2+\delta^2)}{a}\Bigl(\|R\|_{op}+L_sc_1\Bigr)+\delta
            \leq c_3 (\eps+\delta),
        \end{split}
    \end{align}
    where $c_3=\frac{2\sqrt{2}\kappa C'(\eps+\delta)}{a}(\|R\|_{op}+L_sc_1)+1.$ For the third term, 
    \begin{align*}
        \begin{split}
             &\E V\biggl( \bigl(I-\mathscr{K}(\widehat{\Sigma}_j^s)H \bigr) \Bigl(\Psi(P_{\B_V}\widehat{m}_{j-1}^s)-\Psi(P_{\B_V}m^s_{j-1})\Bigr)\biggr)\\
             & \leq \E V\biggl( \bigl(P-\mathscr{K}(\widehat{\Sigma}_j^s)H \bigr) \Bigl(\Psi(P_{\B_V}\widehat{m}_{j-1}^s)-\Psi(P_{\B_V}m^s_{j-1})\Bigr)\biggr)\\
             & +\E V\biggl( (I-P) \Bigl(\Psi(P_{\B_V}\widehat{m}_{j-1}^s)-\Psi(P_{\B_V}m^s_{j-1})\Bigr)\biggr)\\
             & \leq \E \biggl[\mathbb{E}^\mathsf{U_{j-1}^s}\Bigl[\frac{\eps^2+\delta^2}{\lambda_{\min}(H\widehat{\Sigma}_j^sH^*)}\Bigl(\|R\|_{op}+L_sc_1\Bigr) \Bigr]\mathbb{E}^\mathsf{U_{j-1}^s} \Bigl[\Psi(P_{\B_V}\widehat{m}^s_{j-1})-\Psi(P_{\B_V}m^s_{j-1})\Bigr]\biggr]\\
             & + \alpha^{1/2}\E V\bigl(\widehat{m}^s_{j-1}-m^s_{j-1}\bigr)\\
             & \leq \frac{2C'(\eps^2+\delta^2)}{a}\Bigl(\|R\|_{op}+L_sc_1\Bigr)L^{1/2}\E V\bigl(\widehat{m}^s_{j-1}-m^s_{j-1}\bigr)+\alpha^{1/2}\E V\bigl(\widehat{m}^s_{j-1}-m^s_{j-1}\bigr)\\
             & = \Bigl(\frac{c_4}{a}+\alpha^{1/2}\Bigr)\E V\bigl(\widehat{m}^s_{j-1}-m^s_{j-1}\bigr),
        \end{split}
    \end{align*}
    with $c_4=2C'(\eps^2+\delta^2)\Bigl(\|R\|_{op}+L_sc_1\Bigr)L^{1/2}$. Taking $a\geq \frac{2c_4}{1-\alpha^{1/2}}$ yields
    \begin{align}\label{eq:fp surrogate bound 4}
        \E V\biggl( \bigl(I-\mathscr{K}(\widehat{\Sigma}_j^s)H \bigr) \Bigl(\Psi(P_{\B_V}\widehat{m}_{j-1}^s)-\Psi(P_{\B_V}m^s_{j-1})\Bigr)\biggr) \leq \alpha_* \E V\bigl(\widehat{m}^s_{j-1}-m^s_{j-1}\bigr),
    \end{align}
    with $\alpha_*:=\frac{1+\alpha^{1/2}}{2}\in (0,1).$ For the fourth term, we have

\begin{align}\label{eq:fp surrogate bound 5}
        \begin{split}
            &\E V\biggl( \bigl(I-\mathscr{K}(\widehat{\Sigma}_j^s)H \bigr) \Bigl(\Psi(P_{\B_V}m^s_{j-1})-\Psi^s(P_{\B_V}m^s_{j-1}) \Bigr)\biggr)\\
            &\leq \E V\biggl(\bigl(P-\mathscr{K}(\widehat{\Sigma}_j^s)H \bigr) \Bigl(P\Psi(P_{\B_V}m^s_{j-1})-P\Psi^s(P_{\B_V}m^s_{j-1}) \Bigr) \biggr)\\
            &+\E V\biggl((I-P)\Psi(P_{\B_V}m^s_{j-1})-(I-P)\Psi^s(P_{\B_V}m^s_{j-1})\biggr)\\
            & \leq \E \biggl[\mathbb{E}^\mathsf{U_{j-1}^s} \Bigl[\frac{\eps^2+\delta^2}{\lambda_{\min}(H\widehat{\Sigma}_j^sH^*)}\Bigl(\|R\|_{op}+L_sc_1\Bigr) \Bigr]\mathbb{E}^\mathsf{U_{j-1}^s}\Bigl[\Psi(P_{\B_V}m^s_{j-1})-\Psi^s(P_{\B_V}m^s_{j-1})\Bigr]\biggr]\\
            &+ \E V\Bigl((I-P)\Psi(P_{\B_V}m^s_{j-1})-(I-P)\Psi^s(P_{\B_V}m^s_{j-1})\Bigr)\\
            &\leq \E \biggl[\mathbb{E}^\mathsf{U_{j-1}^s} \Bigl[\frac{\eps^2+\delta^2}{\lambda_{\min}(H\widehat{\Sigma}_j^sH^*)}\Bigl(\|R\|_{op}+L_sc_1\Bigr) \Bigr]\sqrt{2}\kappa \biggr]+\delta\\
            & \leq \frac{2\sqrt{2}\kappa C'(\eps^2+\delta^2)}{a}\Bigl(\|R\|_{op}+L_sc_1\Bigr)+\delta
            \leq c_3 (\eps+\delta).
        \end{split}
    \end{align}
    
    Then, we bound the fifth term as
    \begin{align}\label{eq:fp surrogate bound 6}
        \begin{split}
            &\E V\biggl( \bigl(I-\mathscr{K}(\widehat{\Sigma}_j^s)H \bigr) \Bigl(\Psi^s(P_{\B_V}m^s_{j-1})-\mu^s_{j} \Bigr)\biggr)\\
            & \leq \E \biggl[\mathbb{E}^\mathsf{U_{j-1}^s} \Bigl[1+\frac{\eps^2+\delta^2}{\lambda_{\min}(H\widehat{\Sigma}_j^sH^*)}\Bigl(\|R\|_{op}+L_sc_1\Bigr) \Bigr]\mathbb{E}^\mathsf{U_{j-1}^s}\Bigl[V\Bigl(\Psi^s(P_{\B_V}m^s_{j-1})-\mu^s_j\Bigr)\Bigr]\biggr]\\
            & \leq \Bigl(1+\frac{2C'(\eps^2+\delta^2)}{a}(\|R\|_{op}+L_sc_1)\Bigr)\E \biggl[V\Bigl(\Psi^s(P_{\B_V}m^s_{j-1})-\mu^s_j\Bigr)\biggr]\\
            & \leq \Bigl(1+\frac{2C'(\eps^2+\delta^2)}{a}(\|R\|_{op}+L_sc_1)\Bigr)L_s^{1/2}\E \Bigl[\bigl(\Tr(C^s_{j-1})+\Tr(HC^s_{j-1}H^*)\bigr)^{1/2}\Bigr]\\
            & \leq \Bigl(1+\frac{2C'(\eps^2+\delta^2)}{a}(\|R\|_{op}+L_sc_1)\Bigr)L_s^{1/2}\sqrt{2c_1}\eps
             \leq c_5(\eps+\delta),
        \end{split}
    \end{align}
    where $c_5=\Bigl(1+\frac{2C'(\eps^2+\delta^2)}{a}(\|R\|_{op}+L_sc_2)\Bigr)L_s^{1/2}\sqrt{2c_1}.$ For the last term, we have that
    \begin{align*}
        \E V\biggl( \bigl(I-\mathscr{K}(\widehat{\Sigma}_j^s)H \bigr)\frac{1}{N} \sum_{n=1}^N \xi^{(n)}\biggr) & = \E V\biggl( \bigl(P-\mathscr{K}(\widehat{\Sigma}_j^s)H \bigr)\frac{1}{N} \sum_{n=1}^N \xi^{(n)}\biggr)\\
        & \leq (\|R\|_{op}+L_sc_1)\E \biggl[ \frac{\eps^2+\delta^2}{\lambda_{\min}(H\widehat{\Sigma}_j^sH^*)}V \biggl(\frac{1}{N} \sum_{n=1}^N \xi^{(n)}\biggr)\biggr].
    \end{align*} 
    Applying the Cauchy-Schwarz inequality and Lemma \ref{lemma:expected inverse eigenvalue}, it follows that
    \begin{align*}
        &\E V\biggl( \bigl(I-\mathscr{K}(\widehat{\Sigma}_j^s)H \bigr)\frac{1}{N} \sum_{n=1}^N \xi^{(n)}\biggr)\\
        & \leq \bigl(\|R\|_{op}+L_sc_1\bigr)(\eps^2+\delta^2)\E\biggl[\frac{1}{\lambda_{\min}(\widehat{H\Sigma}_j^sH^*)^2}\biggr]^{1/2}\E \biggl[V^2\bigl(\frac{1}{N}\sum_{n=1}^N\xi^{(n)}_j\bigr)\biggr]^{1/2}\\
            & \leq \bigl(\|R\|_{op}+L_sc_1\bigr)(\eps^2+\delta^2)\frac{2C'}{a}\E \biggl[V^2\bigl(\frac{1}{N}\sum_{n=1}^N\xi^{(n)}_j\bigr)\biggr]^{1/2}.
    \end{align*}
    Using the fact that $\E V^2\bigl(\frac{1}{N}\sum_{n=1}^N\xi^{(n)}_j\bigr)= 2 \Tr(\frac{a}{N}I_{k})$ and the assumption that $N\geq 6k$, we have that
            \begin{align} \label{eq:fp surrogate bound 7}
            \begin{split}
            \E V\biggl( \bigl(I-\mathscr{K}(\widehat{\Sigma}_j^s)H \bigr)\frac{1}{N} \sum_{n=1}^N \xi^{(n)}\biggr) & \leq \bigl(\|R\|_{op}+L_sc_1\bigr)(\eps^2+\delta^2)\frac{2C'}{a} \bigl(\frac{2ak}{N}\bigr)^{1/2}\\
             & \leq \bigl(\|R\|_{op}+L_sc_1\bigr)(\eps^2+\delta^2)\frac{2C'}{a^{1/2}} +
              \leq c_6(\eps+\delta),
        \end{split}
    \end{align}
    where we define $c_6:=\bigl(\|R\|_{op}+L_sc_1\bigr)(\eps+\delta)\frac{2C'}{a^{1/2}}$.
    It remains to bound the second term in the right-hand side of \eqref{eq:fp surrogate bound 1}:
    \begin{align*}
        \begin{split}
            &\E V\Bigl(\bigl(I -\mathscr{K}(\widehat{\Sigma}_j^s)H-I +\mathscr{K}(\Sigma^s_j)H \bigr) \bigl(\mu^s_j-H^*y_j \bigr)\Bigr)\\
            & \leq\E V\Bigl(\bigl(I -\mathscr{K}(\widehat{\Sigma}_j^s)H-I +\mathscr{K}(\Sigma^s_j)H \bigr) \bigl(\mu^s_j-\Psi^s(P_{\B_V}m^s_{j-1})\bigr)\Bigr)\\
            & +\E V\Bigl(\bigl(I -\mathscr{K}(\widehat{\Sigma}_j^s)H-I +\mathscr{K}(\Sigma^s_j)H \bigr) \bigl(\Psi^s(P_{\B_V}m^s_{j-1})-H^*y_{j}\bigr)\Bigr).
        \end{split}
    \end{align*}
    We now bound each of these terms. Choosing $a\geq \frac{10NL_sc_2(\eps^2+\delta^2)}{k}$, we have 

    \begin{align}\label{eq:fp surrogate bound 8}
        \begin{split}
            &\E V\Bigl(\bigl(I -\mathscr{K}(\widehat{\Sigma}_j^s)H-I +\mathscr{K}(\Sigma^s_j)H \bigr) \bigl(\mu^s_j-\Psi^s(P_{\B_V}m^s_{j-1})\bigr)\Bigr)\\
            & = \E V\Bigl(\bigl(P -\mathscr{K}(\widehat{\Sigma}_j^s)H-P +\mathscr{K}(\Sigma^s_j)H \bigr) \bigl(\mu^s_j-\Psi^s(m^s_{j-1})\bigr)\Bigr)\\
            & \leq \E \bigl[\bigl(\|P-\mathscr{K}(\widehat{\Sigma}_j^s)H\|_{V,V}+\|P-(\Sigma_j^s)H\|_{V,V}\bigr)V\bigl(\mu^s_j-\Psi^s(P_{\B_V}m^s_{j-1})\bigr) \bigr]\\
            & \leq \E \bigl[\bigl(\|P-\mathscr{K}(\widehat{\Sigma}_j^s)H\|_{V,V}+\|P-(\Sigma_j^s)H\|_{V,V}\bigr)L_s^{1/2}\bigl(\Tr(C_{j-1}^s)+\Tr(HC_{j-1}^sH^*)\bigr)^{1/2} \bigr]\\
            & \leq \E\bigl[\bigl(\|P-\mathscr{K}(\widehat{\Sigma}_j^s)H\|_{V,V}+\|P-(\Sigma_j^s)H\|_{V,V}\bigr)\bigr]L_s^{1/2}\sqrt{2c_1}(\eps+\delta)\\
            & \leq (\|R\|_{op}+L_sc_1)\Bigl(\E\bigl[\frac{\eps^2+\delta^2}{\lambda_{\min}(H\widehat{\Sigma}_j^sH^*)}\bigr] +  \frac{\eps^2+\delta^2}{a}\Bigr)L_s^{1/2}\sqrt{2c_1}(\eps+\delta)\\
            & \overset{\text{(i)}}{\leq} (\|R\|_{op}+L_sc_1)(2C'+1)\frac{\eps^2+\delta^2}{a}L_s^{1/2}\sqrt{2c_1}(\eps+\delta)\leq c_7(\eps+\delta),\\
        \end{split}
    \end{align}
    where (i) uses Lemma \ref{lemma:expected inverse eigenvalue} and $c_7:= (\|R\|_{op}+L_sc_1)(2C'+1)\frac{\eps^2+\delta^2}{a}L_s^{1/2}\sqrt{2c_1}.$ To bound the second term, the triangle inequality and \ref{eq:mf P minus Kalman gain surrogate} yield
    \begin{align*}
        &\E V\Bigl(\bigl(I -\mathscr{K}(\widehat{\Sigma}_j^s)H-I +\mathscr{K}(\Sigma^s_j)H \bigr) \bigl(\Psi^s(P_{\B_V}m^s_{j-1})-H^*y_{j}\bigr)\Bigr)\\
        & = \E V\Bigl(\bigl(P -\mathscr{K}(\widehat{\Sigma}_j^s)H-P +\mathscr{K}(\Sigma^s_j)H \bigr) \bigl(P\Psi^s(P_{\B_V}m^s_{j-1})-P\Psi(u_{j-1})+H^*\eta_j\bigr)\Bigr)\\
        & \leq \E \bigl[\bigl(\|P-\mathscr{K}(\widehat{\Sigma}^s_j)\|_{V,V}+\|P-\mathscr{K}(\Sigma^s_j)\|_{V,V}\bigr)V\bigl(P\Psi^s(P_{\B_V}m^s_{j-1})-P\Psi(u_{j-1})+H^*\eta_j\bigr)\bigr]\\
        & \leq \E \bigl[\bigl(\|P-\mathscr{K}(\widehat{\Sigma}^s_j)\|_{V,V}+\|P-\mathscr{K}(\Sigma^s_j)\|_{V,V}\bigr)V\bigl(P\Psi^s(P_{\B_V}m^s_{j-1})-P\Psi(P_{\B_V}m^s_{j-1})\bigr)\bigr]\\
        &+ \E \bigl[\bigl(\|P-\mathscr{K}(\widehat{\Sigma}^s_j)\|_{V,V}+\|P-\mathscr{K}(\Sigma^s_j)\|_{V,V}\bigr)V\bigl(P\Psi(P_{\B_V}m^s_{j-1})-P\Psi(u_{j-1})+H^*\eta_j\bigr)\bigr]\\
        &\leq (\|R\|_{op}+L_sc_1)\E \bigl[\bigl(\frac{\eps^2+\delta^2}{\lambda_{\min}(H\widehat{\Sigma}^s_jH^*)}+\frac{\eps^2+\delta^2}{a}\bigr)V\bigl(P\Psi^s(P_{\B_V}m^s_{j-1})-P\Psi(P_{\B_V}m^s_{j-1})\bigr)\bigr]\\
        & + (\|R\|_{op}+L_sc_1)\E \bigl[\bigl(\frac{\eps^2+\delta^2}{\lambda_{\min}(H\widehat{\Sigma}^s_jH^*)}+\frac{\eps^2+\delta^2}{a}\bigr)V\bigl(P\Psi(P_{\B_V}m^s_{j-1})-P\Psi(u_{j-1})+H^*\eta_j\bigr)\bigr]\\
        & \leq (\|R\|_{op}+L_sc_1)\E \bigl[\bigl(\frac{\eps^2+\delta^2}{\lambda_{\min}(H\widehat{\Sigma}^s_jH^*)}+\frac{\eps^2+\delta^2}{a}\bigr)\kappa \bigr]\\
        & + (\|R\|_{op}+L_sc_1)\E \bigl[\bigl(\frac{\eps^2+\delta^2}{\lambda_{\min}(H\widehat{\Sigma}^s_jH^*)}+\frac{\eps^2+\delta^2}{a}\bigr)V\bigl(\Psi(P_{\B_V}m^s_{j-1})-\Psi(u_{j-1})+H^*\eta_j\bigr)\bigr],
    \end{align*}
    where the last inequality uses Assumption \ref{assumption:surrogate model assumption}.
    Next, by conditional independence and Lemma \ref{lemma:expected inverse eigenvalue},
    \begin{align*}
        &\E V\Bigl(\bigl(I -\mathscr{K}(\widehat{\Sigma}_j^s)H-I +\mathscr{K}(\Sigma^s_j)H \bigr) \bigl(\Psi^s(P_{\B_V}m^s_{j-1})-H^*y_{j}\bigr)\Bigr)\\
        & \leq (\|R\|_{op}+L_sc_1)\bigl(\E \bigl[\frac{\eps^2+\delta^2}{\lambda_{\min}(H\widehat{\Sigma}^s_jH^*)}\bigr]+\frac{\eps^2+\delta^2}{a}\bigr)\kappa\\
        & + (\|R\|_{op}+L_sc_1)\E \Bigl[\Bigl(\E^{\mathsf{U}_{j-1}^s}\bigl[\frac{\eps^2+\delta^2}{\lambda_{\min}(H\widehat{\Sigma}^s_jH^*)}\bigr]+\frac{\eps^2+\delta^2}{a}\Bigr)\E^{\mathsf{U}_{j-1}^s}\bigl[V\bigl(\Psi(P_{\B_V}m^s_{j-1})-\Psi(u_{j-1})\bigr)\bigr] \Bigr]\\
        & + (\|R\|_{op}+L_sc_1)\Bigl(\E\bigl[\frac{\eps^2+\delta^2}{\lambda_{\min}(H\widehat{\Sigma}^s_jH^*)}\bigr]+\frac{\eps^2+\delta^2}{a}\Bigr)\E V\bigl(H^*\eta_j\bigr)\\
        & \leq (\|R\|_{op}+L_sc_1)(2C'+1)\frac{\eps^2+\delta^2}{a}\kappa\\
        & + (\|R\|_{op}+L_sc_1)(2C'+1)\frac{\eps^2+\delta^2}{a}\E\bigl[V\bigl(m^s_{j-1}-u_{j-1}\bigr)\bigr]\\
        & + (\|R\|_{op}+L_sc_1)(2C'+1)\frac{\eps^2+\delta^2}{a}\sqrt{2}\Tr(R)^{1/2}\eps
    \end{align*}
    where the last inequality also uses the fact that $\Psi$ is locally Lipschitz on $\B_V$ and that $\E V(H^*\eta_j)\leq \sqrt{2}\Tr(R)^{1/2}\eps$. Since $j>j^*$, we have that $\E V(m^s_{j-1}-u_{j-1})\leq 2(\mathsf{C}_1+\mathsf{C}_3)(\eps+\delta),$ and thus
    \begin{align}\label{eq:fp surrogate bound 9}
        \begin{split}
             &\E V\Bigl(\bigl(I -\mathscr{K}(\widehat{\Sigma}_j^s)H-I +\mathscr{K}(\Sigma^s_j)H \bigr) \bigl(\Psi^s(P_{\B_V}m^s_{j-1})-H^*y_{j}\bigr)\Bigr)\\
             & \leq (\|R\|_{op}+L_sc_1)(2C'+1)\bigl(\kappa\frac{\eps+\delta}{a}+2\frac{(\eps+\delta)^2}{a}(\mathsf{C}_1+\mathsf{C}_3)+\frac{\eps+\delta}{a}\sqrt{2}\Tr(R)^{1/2}\eps\bigr)(\eps+\delta)\\
             & =c_8(\eps+\delta),
        \end{split}
    \end{align}
    where we define $c_8=(\|R\|_{op}+L_sc_1)(2C'+1)\bigl(\kappa\frac{\eps+\delta}{a}+2\frac{(\eps+\delta)^2}{a}(\mathsf{C}_1+\mathsf{C}_3)+\frac{\eps+\delta}{a}\sqrt{2}\Tr(R)^{1/2}\eps\bigr).$
   Combining \eqref{eq:fp surrogate bound 2}, \eqref{eq:fp surrogate bound 3}, \eqref{eq:fp surrogate bound 4}, \eqref{eq:fp surrogate bound 5}, \eqref{eq:fp surrogate bound 6}, \eqref{eq:fp surrogate bound 7}, \eqref{eq:fp surrogate bound 8}, and \eqref{eq:fp surrogate bound 9} we have
    \begin{align}
        \E V\bigl(\widehat{m}_j^s-m^s_j\bigr) \leq \alpha_*\E V\bigl(\widehat{m}_{j-1}^s-m^s_{j-1}\bigr) +(c_2+2c_3+c_5+c_6+c_7+c_8)(\eps+\delta).
    \end{align}
    The discrete Gronwall inequality and the equivalence of $V(\cdot)$ and $\|\cdot\|$ yield the desired result.
    \end{proof}
\end{theorem}

\section{Numerical experiments}\label{sec:numerics}
In this section, we numerically verify our theoretical results with the Lorenz-96 system given in Example \ref{ex:Lorenz96}. Throughout our experiments, we take the state dimension to be $d=60$ and the observation time-step to be $\Delta t=0.1$. As described in Example \ref{ex:Lorenz96}, we take our observation operator to be the identity matrix with every third row removed, and thus $k=40.$

\begin{figure}
    \centering
    \includegraphics[width=.9\linewidth]{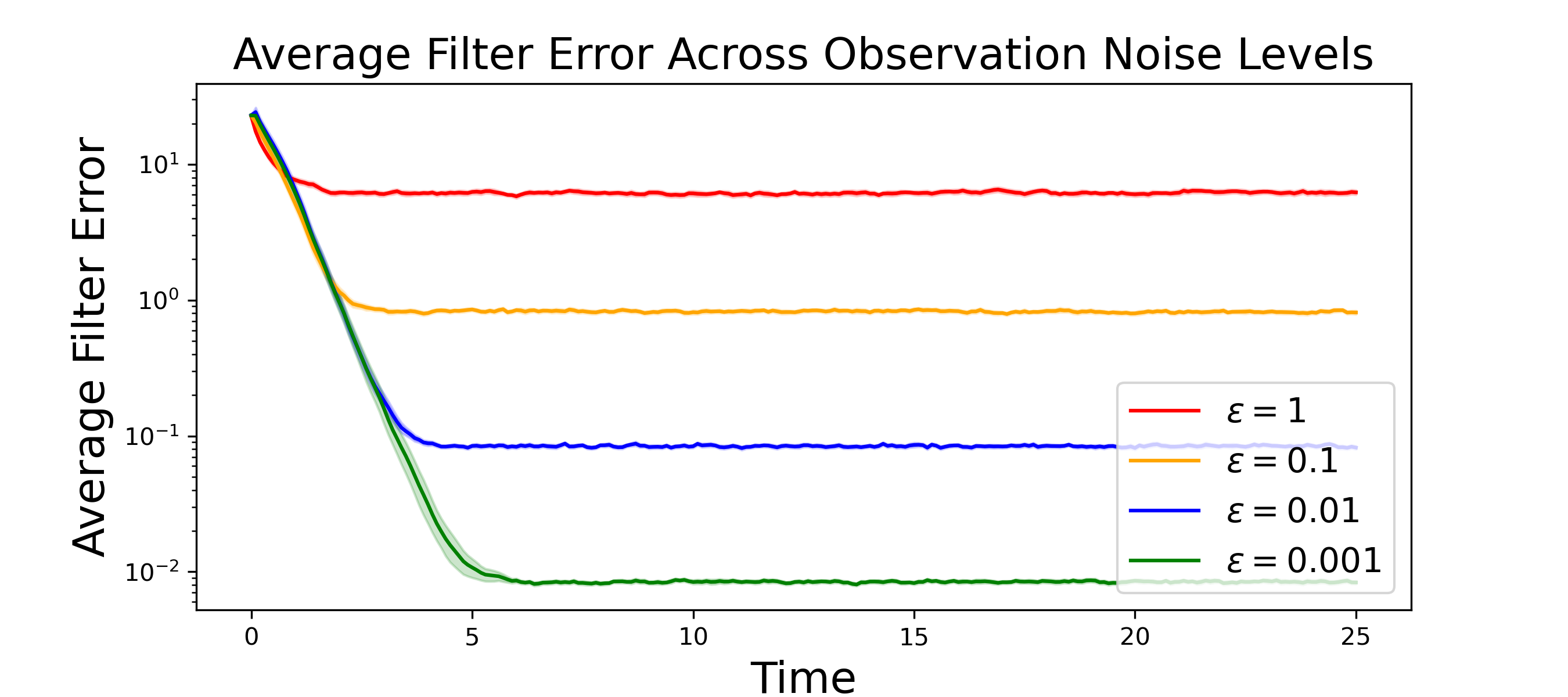}
    \caption{Average filter error over 50 Monte Carlo trials for decreasing observation noise levels, plotted with two standard errors (shaded).} 
    \label{fig:decreasing noise plot}
\end{figure}

    To illustrate the filter accuracy result given in Theorem \ref{th:main1}, we compare the performance of Algorithm \ref{algorithm:finite ensemble} for decreasing noise levels. For 50 Monte Carlo trials, we initialize the true signal at a random Gaussian vector $u_0\sim \Nc(0,10 I)$ and evolve the true state until time $T=25$. We generate observations $y_j$ at noise levels $\eps=(1,10^{-1},10^{-2},10^{-3}).$ Taking $N=50$ ensemble members, we initialize the ensemble as $\{u_0^{(n)}\}_{n=1}^N\sim \Nc(0,I)$. We then run the Ensemble Transform Kalman Filter with variance inflation parameter $a=1$. In Figure \ref{fig:decreasing noise plot}, we plot the average error between the analysis mean and the true signal, $\|\widehat{m}_j-u_j\|$, over the 50 Monte Carlo trials. As predicted by our theory, the errors are proportional to the noise level, $\epsilon.$

To illustrate Theorem \ref{th:main2}, we compare the performance of Algorithm \ref{algorithm:finite ensemble surrogate} with surrogate models of increasing fidelity. To generate training data for our surrogate models, we sample $10^3$ standard normal random vectors and integrate the Lorenz-96 system until time $T=100.$ This generates $10^6$ input-output pairs of the dynamics map that we take to be our training data set. We consider a convolutional neural network with architecture similar to that used in \cite{brajard2020combining} with $30,003$ total trainable parameters. The structure of the neural network is detailed in Figure \ref{fig:network details figure}. The neural network weights are initialized randomly using Pytorch's default initialization. To train the neural network, we use the mean squared error loss function with the Adam \cite{kingma2014adam} optimizer with a learning rate of $0.0003$, momentum $0.1$, and a batch size of $2000.$ We train three neural networks using different amounts of training data, meant to represent a low-fidelity, medium-fidelity, and high-fidelity machine model. The size of the training set and number of epochs for each model are summarized in Table \ref{table:ML summary}.
\begin{figure}
    \centering
    \begin{subfigure}{1\textwidth}

\begin{tikzpicture}[
    font=\small,
    every node/.style={draw, minimum height=1cm, minimum width=1cm, align=center},
    arr/.style={->, thick},
    arr2/.style={-, thick}
]

\node[circle] (u) at (0, 0) {$u$};

\node (CNN11) at (2, 2) {CNN$^{(1)}_1$};
\node (CNN12) at (2, 0) {CNN$^{(2)}_1$};
\node (CNN13) at (2, -2) {CNN$^{(3)}_1$};

\node[circle] (mul) at (4, -1) {$\times$};

\node (CNN2) at (6, 0) {CNN$_2$};

\node (CNN3) at (8, 0) {CNN$_3$};

\node (CNN4) at (10, 0) {CNN$_4$};

\node (CNN5) at (12, 0) {CNN$_5$};

\node[circle] (f) at (14, 0) {$\Psi^s(u)$};

\draw[arr] (u) -- (CNN11);
\draw[arr] (u) -- (CNN12);
\draw[arr] (u) -- (CNN13);

\draw[arr] (CNN11) -- (CNN2);
\draw[arr] (CNN12) -- (mul);
\draw[arr] (CNN13) -- (mul);

\draw[arr] (mul) -- (CNN2);
\draw[arr] (CNN2) -- (CNN3);
\draw[arr] (CNN3) -- (CNN4);
\draw[arr] (CNN4) -- (CNN5);
\draw[arr] (CNN5) -- (f);

\end{tikzpicture}
    \caption{Network structure.}
    \label{fig:network structure}
    \end{subfigure}
    \begin{subfigure}{1\textwidth}
	\begin{center}
		\begin{tabular}{ | c | c |c | c |}
			\hline
			 & \#in  & \#out   & kernel size \\ \hline
		$\text{CNN}_1$ & 1 & 72(24$\times$ 3) & 5, circular  \\ \hline
         $\text{CNN}_2$ & 48(24$\times$ 2) & 37 & 5, circular  \\ \hline
          $\text{CNN}_3$ & 37 & 37 & 5, circular  \\ \hline
          $\text{CNN}_4$ & 37 & 37 & 5, circular  \\ \hline
          $\text{CNN}_5$ & 37 & 1 & 1  \\ \hline
		\end{tabular}
	\end{center}
    \caption{Network details.}

\label{fig:network details}
    \end{subfigure}
    \caption{The structure of $\Psi^s$. The output channels of $\text{CNN}_1$ are divided into three groups of equal length, $\text{CNN}_1^{(1)}$,  $\text{CNN}_1^{(2)}$, and  $\text{CNN}_1^{(3)}$. The input channels to $\text{CNN}_2$ are a concatenation of  $\text{CNN}_1^{(1)}$ and $(\text{CNN}_1^{(1)}\times \text{CNN}_1^{(2)})$, where the multiplication is point-wise. }
    \label{fig:network details figure}
\end{figure}
To estimate the model error in the unobserved components of the system, $\delta$, defined in Assumption \ref{assumption:surrogate model assumption}, we select $10^3$ points uniformly at random from the attractor by numerically integrating a long run of the dynamics and compute the maximum error of each of the models over the $10^3$ points. The estimated model errors are reported in Table \ref{table:ML summary}. 

\begin{figure}
    \centering
    \includegraphics[width=1\linewidth]{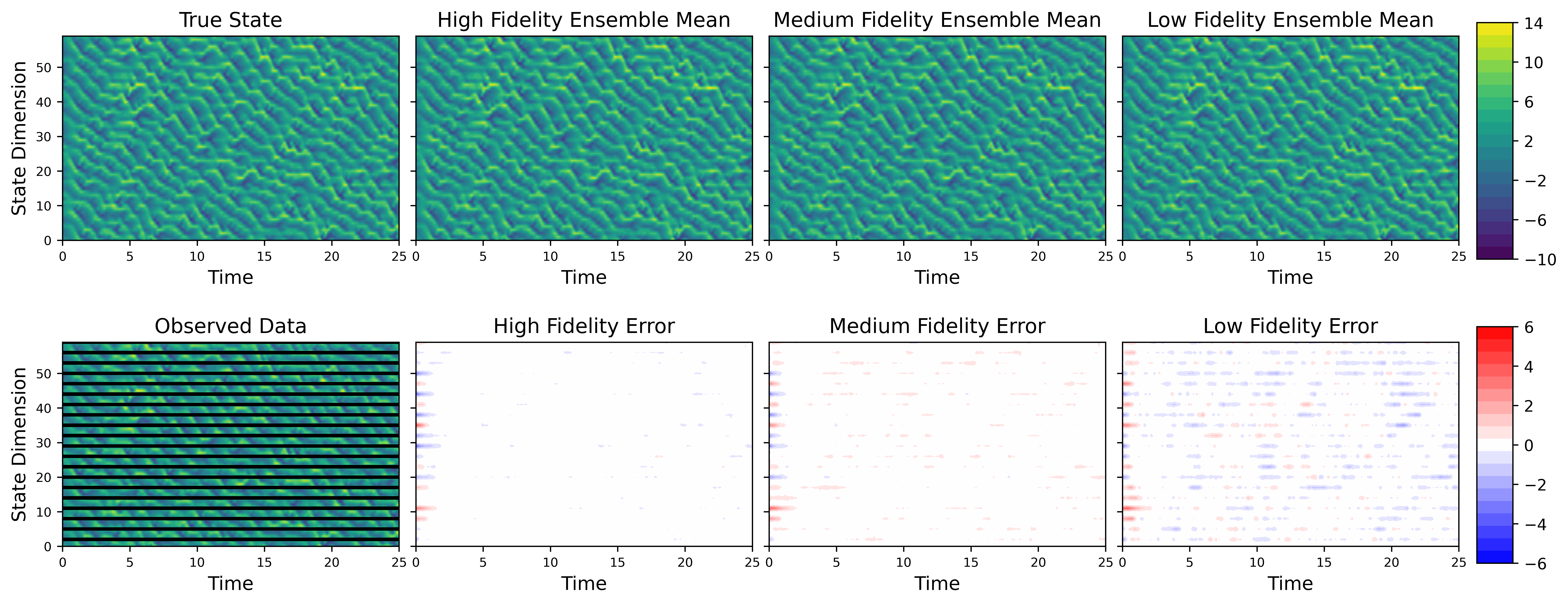}
    \caption{Visualization of the true state trajectory over time, the ensemble means, the observed data, and the difference between the true state and the ensemble means. }
    \label{fig:state plot}
\end{figure}

In Figure \ref{fig:state plot}, we visualize the performance of the Ensemble Transform Kalman filter with surrogate models, implemented with $N=50$ ensemble members and variance inflation parameter $a=10$. For a single realization of the true signal, we plot the evolution of the true state until time $T=25$, along with the observed data, the surrogate model filter means, and the difference between the true state and the filter estimates. Visually, the ensemble means output by each of the filters appear to track the true signal quite well. The error plots reveal that, as predicted by our theory, filtering with the low fidelity surrogate model incurs the largest error, and the high fidelity error incurs the smallest error.

\begin{figure}
    \centering
    \includegraphics[width=1\linewidth]{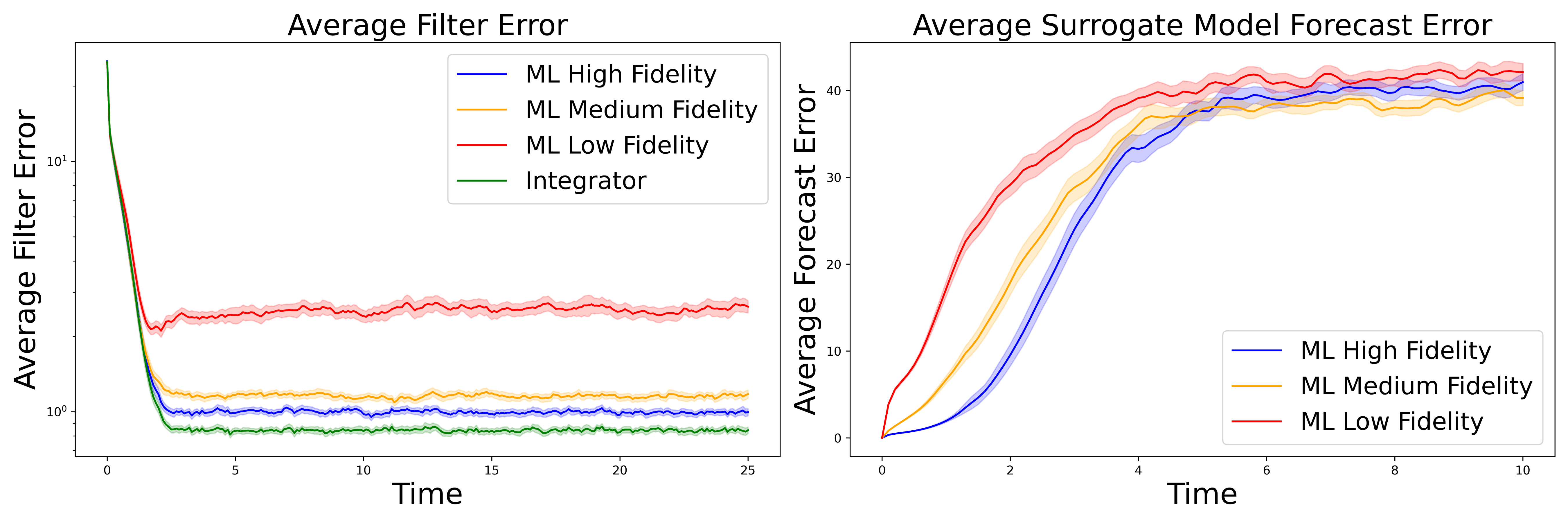}
    \caption{The left plot displays the average filter error over 50 Monte Carlo trials plotted with two standard errors. The right plot displays the average surrogate model forecast error from an exactly known initial condition over 50 Monte Carlo trials plotted with two standard errors. This figure shows that long-time filter accuracy is possible with surrogate models that only give accurate short-term forecasts.}
    \label{fig:surrogate model experiments}
\end{figure}

Next, for 50 Monte Carlo trials, we initialize the true signal at a random Gaussian vector $u_0\sim \Nc(0,10 I)$ and we initialize the ensemble as $\{u_0^{(n)}\}_{n=1}^N\sim \Nc(0,I)$. We then run the Ensemble Transform Kalman Filter with $N=50$ ensemble members and variance inflation parameter $a=10$. In the first plot in Figure \ref{fig:surrogate model experiments}, we display the average error in the analysis forecast mean $\widehat{m}_j$ from the true state $u_j$. We observe that the average filter error is proportional to the model error in the unobserved components, as predicted by our theory. The filter errors averaged over times $T=10$ to $T=25$ are reported in Table \ref{table:ML summary}. 

 We then consider the forecast error from the true dynamics incurred by applying the surrogate model sequentially from an exactly known initial condition to emulate the dynamics over moderate time-scales. In the second plot in Figure \ref{fig:surrogate model experiments}, we display the average error over 50 Monte Carlo trials with initial conditions drawn from $u_0\sim \Nc(0,10 I)$. By time $T=4$, the error in all surrogate model predictions is large. We see that despite the surrogate models being unable to accurately forecast the dynamics beyond short timescales, they can still be used within the ensemble Kalman filter to accurately track the signal given sufficiently informative observations.

\begin{table}[H]
	\begin{center}
		\begin{tabular}{ | c | c |c | c | c |}
			\hline
			Fidelity & Training Samples  & Epochs   & Model Error & Filter Error  \\ \hline
		 Low & $10^3$ & 20 & $\delta\approx 2.16$ & 2.59 \\ \hline
         Medium & $10^4$ & 50 & $\delta\approx 1.02 $ & 1.16 \\ \hline
         High & $10^6$ & 200 & $\delta\approx 0.35$ & 0.99 \\ \hline
		\end{tabular}
		\caption{Summary of the training details, estimated model error, and average filter error for each model. }	
        \label{table:ML summary}
	\end{center}
    
\end{table}

\section{Conclusions}\label{sec:conclusions}
This paper established long-time accuracy of ensemble Kalman filters, introducing conditions on the dynamics and the observations under which the estimation error remains small in the long-time horizon. We have shown that these conditions can be verified for many partially-observed dynamical systems, such as the Lorenz models and the 2-dimensional Navier-Stokes equations. In addition, we have rigorously established and empirically demonstrated that ensemble Kalman filters can provide long-time accurate state estimation even when employing machine-learned surrogate models that only give accurate short-term forecasts. 

There are several natural research directions that stem from this work. First, our analysis focuses on linear observations; the nonlinear case could be studied by working on an extended state-space, as in \cite{anderson2001ensemble}. Second, it would be interesting to extend the theory beyond the case of autonomous dynamical systems and fixed observation operators considered here.  Third, our work does not account for covariance localization, which can improve the empirical performance of the algorithm and reduce the required ensemble size. Lastly, there are numerous opportunities to leverage the ideas in this paper to co-learn the state, the dynamics, and the filtering algorithms in data assimilation tasks across many application domains. 

\section*{Acknowledgments}
This work was partially funded by the NSF CAREER award NSF DMS-2237628. The authors are grateful to Melissa Adrian, Omar Al-Ghattas, Mihai Anitescu, and Nisha Chandramoorthy for helpful discussions. 
\bibliographystyle{siamplain}
\bibliography{references}

\end{document}